\documentclass[oneside, 12pt]{amsart}

\usepackage{amsfonts}
\usepackage{amssymb}
\usepackage{latexsym}
\usepackage{graphics}
\usepackage{amsthm}
\usepackage{amsmath}
\usepackage[all]{xy}

\begin{document}

\newcommand{\norm}[1]{\| #1 \|}
\def\N{\mathbb N}
\def\Z{\mathbb Z}
\def\Q{\mathbb Q}
\def\mod{\textit{\emph{~mod~}}}
\def\R{\mathcal R}
\def\S{\mathcal S}
\def\*  C{{*  \mathcal C}}
\def\C{\mathcal C}
\def\D{\mathcal D}
\def\J{\mathcal J}
\def\M{\mathcal M}
\def\T{\mathcal T}

\newcommand{\Hom}{{\rm Hom}}
\newcommand{\End}{{\rm End}}
\newcommand{\Ext}{{\rm Ext}}
\newcommand{\Mor}{{\rm Mor}\,}
\newcommand{\Aut}{{\rm Aut}\,}
\newcommand{\Hopf}{{\rm Hopf}\,}
\newcommand{\Ann}{{\rm Ann}\,}
\newcommand{\Ker}{{\rm Ker}}
\newcommand{\Coker}{{\rm Coker}\,}
\newcommand{\im}{{\rm Im}\,}
\newcommand{\coim}{{\rm Coim}\,}
\newcommand{\Trace}{{\rm Trace}\,}
\newcommand{\Char}{{\rm Char}\,}
\newcommand{\Mod}{{\rm mod}}
\newcommand{\Spec}{{\rm Spec}\,}
\newcommand{\sgn}{{\rm sgn}\,}
\newcommand{\Id}{{\rm Id}\,}
\newcommand{\Com}{{\rm Com}\,}
\newcommand{\codim}{{\rm codim}}
\newcommand{\Mat}{{\rm Mat}}
\newcommand{\can}{{\rm can}}
\newcommand{\sign}{{\rm sign}}
\newcommand{\kar}{{\rm kar}}
\newcommand{\rad}{{\rm rad}}
\newcommand{\dx}{\, \mbox{\rm d}}
\newcommand{\Rad}{{\rm Rad}}

\def\lan{\langle}
\def\ran{\rangle}
\def\ot{\otimes}

\def\id{{\small \textit{\emph{1}}}\!\!1}
\def\To{{\multimap\!\to}}
\def\bigperp{{\LARGE\textrm{$\perp$}}}
\newcommand{\QED}{\hspace{\stretch{1}}
\makebox[0mm][r]{$\Box$}\\}

\def\RR{{\mathbb R}}
\def\FF{{\mathbb F}}
\def\NN{{\mathbb N}}
\def\CC{{\mathbb C}}
\def\DD{{\mathbb D}}
\def\ZZ{{\mathbb Z}}
\def\QQ{{\mathbb Q}}
\def\HH{{\mathbb H}}
\def\units{{\mathbb G}_m}
\def\GG{{\mathbb G}}
\def\EE{{\mathbb E}}
\def\FF{{\mathbb F}}
\def\rightact{\hbox{$\leftharpoonup$}}
\def\leftact{\hbox{$\rightharpoonup$}}

\newcommand{\Aa}{\mathcal{A}}
\newcommand{\Bb}{\mathcal{B}}
\newcommand{\Cc}{\mathcal{C}}
\newcommand{\Dd}{\mathcal{D}}
\newcommand{\Ee}{\mathcal{E}}
\newcommand{\Ff}{\mathcal{F}}
\newcommand{\Hh}{\mathcal{H}}
\newcommand{\Ii}{\mathcal{I}}
\newcommand{\Mm}{\mathcal{M}}
\newcommand{\Pp}{\mathcal{P}}
\newcommand{\Rr}{\mathcal{R}}
\newcommand{\squig}{\rightsquigarrow}

\def\*  C{{}*  \hspace*  {-1pt}{\Cc}}

\def\text#1{{\rm {\rm #1}}}

\def\smashco{\mathrel>\joinrel\mathrel\triangleleft}
\def\cosmash{\mathrel\triangleright\joinrel\mathrel<}

\def\Nat{\dul{\rm Nat}}

\newtheorem{prop}{Proposition}[section]
\newtheorem{lemma}[prop]{Lemma}
\newtheorem{cor}[prop]{Corollary}
\newtheorem{theo}[prop]{Theorem}

\theoremstyle{definition}
\newtheorem{Def}[prop]{Definition}
\newtheorem{ex}[prop]{Example}
\newtheorem{exs}[prop]{Examples}
\newtheorem{Not}[prop]{Notation}
\newtheorem{Ax}[prop]{Axiom}
\newtheorem{rems}[prop]{Remarks}
\newtheorem{rem}[prop]{Remark}

\def\smashco{\mathrel>\joinrel\mathrel\triangleleft}

\title{Measures, states and de Finetti maps on pseudo-BCK algebras}

\begin{abstract}
In this paper, we extend the notions of states and measures
presented in \cite{DvPu} to the case of pseudo-BCK algebras and
study similar properties. We prove that, under some conditions, the
notion of a state in the sense of \cite{DvPu} coincides with the
Bosbach state, and we extend to the case of pseudo-BCK algebras some
results proved by J. K\"uhr only for  pseudo-BCK semilattices. We
characterize extremal states, and show that the quotient pseudo-BCK
algebra over the kernel of a measure can be embedded into the
negative cone of an archimedean $\ell$-group. Additionally, we
introduce a Borel state and using results by J. K\"uhr and D.
Mundici from \cite{Kumu}, we prove a relationship between de Finetti
maps, Bosbach states and Borel states.

\end{abstract}

\author[Ciungu, Dvure\v censkij]{Lavinia Corina Ciungu$^1$, Anatolij Dvure\v censkij$^2$}
\date{}
\thanks{The first author thanks SAIA for the fellowship in Slovakia, Summer
2008, 2009, the second author is thankful for the support by Center
of Excellence SAS -~Quantum Technologies~-, Center of Excellence
QUTE, the grant VEGA No. 2/0032/09 SAV and by the Slovak Research
and Development Agency under the contract No. APVV-0071-06,
Bratislava.}

\maketitle

\section{Introduction}

BCK algebras were introduced originally by K. Is\`eki in \cite{Ise}
with a  binary operation $*$ modeling the set-theoretical difference
and with a constant element $0$ that is a least element. Another
motivation is from classical and non-classical propositional calculi
modeling logical implications. Such algebras contain as a special
subfamily a family of MV-algebras where some important fuzzy
structures can be studied. For more about BCK algebras, see
\cite{Men1}.

Pseudo-BCK algebras were originally introduced by G. Georgescu and
A. Iorgulescu in \cite{Geo15}  as  algebras  with ``two
differences", a left- and right-difference, instead of one $*$ and
with a constant element $0$ as the least element.  In \cite{DvVe}, a
special subclass of pseudo-BCK algebras, called \L ukasiewicz
pseudo-BCK algebras, was introduced and it was shown that it is
always a subalgebra of the positive cone of some $\ell$-group (not
necessarily abelian). The class of \L ukasiewicz pseudo-BCK algebras
is a variety whereas  the class of pseudo-BCK algebras is not; it is
only a quasivariety because it is not closed under homomorphic
images. Nowadays pseudo-BCK algebras are used in a dual form, with
two implications, $\to$ and $\rightsquigarrow$ and with  one
constant element $1$ that is a greatest element. Thus such
pseudo-BCK algebras are in a ``negative cone" and are also called
``left-ones". For a guide through the pseudo-BCK algebras realm, see
the monograph \cite{Ior3}.

States or measures give a probabilistic interpretation of randomness
of events of given algebraic structures. For MV-algebras, Mundici
introduced states (an analogue of probability measures)  in 1995,
\cite{Mun1}, as averaging of the truth-value in \L ukasiewicz logic.
Measures on BCK algebras were introduced by A. Dvure\v censkij in
\cite{Dvu1, DvPu}. Based on the notion of a measure it was defined a
concept of state on these structures if such BCK algebra admits a
smallest element $0.$ Today the notion of states has many forms: The
notion of a Bosbach state has been studied for other algebras of
fuzzy structures such as pseudo-BL algebras, \cite{Geo1}, bounded
non-commutative $R\ell$-monoids, \cite{DvRa1, DvRa2, DvRa3},
residuated lattices, \cite{Ciu1}, pseudo-BCK semilattices and
pseudo-BCK algebras, \cite{Kuhr}.

In the present paper, we study Bosbach states and measures on pseudo
BCK algebras.  In general, it can happen that on bounded BCK
algebras  states  can fail. We will study states, extremal states,
Bosbach states, state-morphisms and we show the relationships among
them. Such relations were studied for pseudo MV-algebras,
\cite{Dvu2}, and generalized by K\"uhr \cite{Kuhr} for pseudo-BCK
algebras that are $\vee$-semilattices. In our paper we show that the
existence of the join in the pseudo-BCK algebra is not substantial
for our study. The main results say that the quotient pseudo-BCK
algebra that is downwards-directed over the kernel of a measure can
be embedded as a subalgebra into the negative cone of an abelian and
archimedean $\ell$-group. In particular that $A$ is with strong
unit, the embedding is even onto. We will apply this result to
characterize extremal state-measures on unital $\ell$-groups.

Finally, we show how Bosbach states and state-measures appear with
respect to de Finetti's coherence principle.

The paper is organized as follows: Section 2 contains preliminary
notions on pseudo-BCK algebras. Section 3 is dedicated to states,
extremal states, state-morphisms and Bosbach states, kernels on
pseudo-BCK algebras. Section 4 deals with a generalization of a
notion of measures, Section 5  presents results on state-measures on
pseudo-BCK algebras with strong unit. The last section deals with de
Finetti's notion of coherence and it gives some integral
representations of states and state-measures via Borel states.

\section{Preliminaries on Pseudo-BCK Algebras}

In the present section, we give elements of theory of pseudo-BCK
algebras.

\begin{Def} \label{psBCK-1.10} $\rm($\cite{Ior1}$\rm)$ A {\it pseudo-BCK algebra}
is a structure ${\mathcal
A}=(A,\leq,\rightarrow,\rightsquigarrow,1)$ where $\leq$ is a
partial binary relation on $A$, $\rightarrow$ and $\rightsquigarrow$
are total binary operations on $A$ and $1$ is an element of $A$
satisfying, for all $x,y,z \in A$, the  axioms:
\begin{enumerate}
\item[$(A_1)$] $x \rightarrow y \leq (y \rightarrow z) \rightsquigarrow
(x \rightarrow z)$, $x \rightsquigarrow y \leq (y \rightsquigarrow
z) \rightarrow (x \rightsquigarrow z);$
\item[$(A_2)$] $x \leq (x \rightarrow y) \rightsquigarrow y$,$\:\:\:$
$x \leq (x \rightsquigarrow y) \rightarrow y;$
\item[$(A_3)$] $x \leq x;$
\item[$(A_4)$] $x \leq 1;$
\item[$(A_5)$] if $x \leq y$ and $y \leq x$, then $x = y;$
\item[$(A_6)$] $x \leq y$ iff $x \rightarrow y = 1$ iff $x \rightsquigarrow
y = 1$.
\end{enumerate}
\end{Def}

Without loss of generality, we will denote a pseudo-BCK algebra
$(A,\leq,\rightarrow,\rightsquigarrow,1)$ simply by $A.$

\begin{rems} \label{psBCK-1.20}  $\rm($\cite{Ior1}$\rm)$
$(1)$ A pseudo-BCK algebra $A$ is a BCK algebra if
$\rightarrow \: = \: \rightsquigarrow$.\\
$(2)$ The partial operation $\leq$ is in fact a partial order on $A$.\\
$(3)$  If there is an element $0$ in $A$ such that $0 \leq x$ (i.e.
$0 \rightarrow x = 0 \rightsquigarrow x=1$), for all $x \in A$, then
$0$ is called the \emph{zero} of $ A$. A pseudo-BCK algebra with
zero is called a \emph{bounded pseudo-BCK algebra} and it is denoted
by ${\mathcal A}=(A,\leq,\rightarrow,\rightsquigarrow,0,1),$ and in
a simple way also as $A.$ In such a case, we define two negations,
$^-$ and $^\sim$, for any element $x \in A:$
$$x^- := x\to 0,\quad x^\sim := x\rightsquigarrow 0.$$
\end{rems}

\begin{ex} \label{psBCK-1.10-20}
Let $(G,+,0,\wedge,\vee)$ be an $\ell$-group (= a lattice ordered
group that
is not necessarily abelian). \\
On the negative cone $G^-=\{g\in G \mid g\leq 0\}$ we define:
\begin{eqnarray*} g\rightarrow h&:=&h-(g\vee h)= (h-g)\wedge 0, \\
 g\rightsquigarrow h&:=&-(g\vee h)+h = (-g+h)\wedge 0.
\end{eqnarray*} Then $(G^-,\rightarrow,\rightsquigarrow,0)$ is a pseudo-BCK
algebra.
\end{ex}

\begin{ex} \label{psBCK-1.10-30} Let $(G,+,0,\wedge,\vee)$ be an
$\ell$-group with a strong unit $u\geq 0$ (i.e. given $g\in G$ there
is an integer $n\ge 1$ such that $g \le nu$).
On the interval $[-u,0]$ we define: \\
$\hspace*{5cm}$ $x\rightarrow y:=(y-x)\wedge 0,$ \\
$\hspace*{5cm}$ $x\rightsquigarrow y:=(-x+y)\wedge 0$. \\
Then $([-u,0],\rightarrow,\rightsquigarrow,-u,0)$ is a bounded
pseudo-BCK algebra with $x^-=-u+x $ and $x^\sim = -x+u.$ In a
similar way, $((-u,0],\to,\rightsquigarrow,0)$ is a pseudo-BCK
algebra that is not bounded.
\end{ex}

\begin{ex} \label{psBCK-1.10-40} Let $(G,+,0,\wedge,\vee)$ be an  $\ell$-group
with a strong unit $u\geq 0$. On the interval $[0,u]$ we define: \\
$\hspace*{5cm}$ $x\rightarrow y:=(u-x+y)\wedge u,$ \\
$\hspace*{5cm}$ $x\rightsquigarrow y:=(y-x+u)\wedge u $. \\
Then $([0,u],\rightarrow,\rightsquigarrow,0,u)$ is a bounded
pseudo-BCK algebra with $x^-=u-x$ and $x^\sim = -x+u.$

If on $[0,u]$ we set $\to_1 =\squig$ and $\squig_1 = \to,$ then
$([0,u],\to_1,\squig_1,0,u)$ is isomorphic with
$[-u,0],\to,\squig,-u,0)$ under the isomorphism $x\mapsto x-u,$ $x
\in [0,u].$
\end{ex}

\begin{prop} \label{psBCK-1.20-5} $\rm($\cite{Ior12}, \cite{Ior13}$\rm)$
In any pseudo-BCK algebra the following properties hold:\\
$(c_1)$ $x \rightarrow (y \rightsquigarrow z) = y \rightsquigarrow (x \rightarrow z)$ and
                $x \rightsquigarrow (y \rightarrow z) = y \rightarrow (x \rightsquigarrow z).$ \\
$(c_2)$ $x \leq y \rightarrow x$, $\:\:$ $x \leq y \rightsquigarrow x$.\\
$(c_3)$ $[(y \rightarrow x) \rightsquigarrow x] \rightarrow x = y \rightarrow x$, $\:\:$
                $[(y \rightsquigarrow x) \rightarrow x] \rightsquigarrow x = y \rightsquigarrow
                x$.\\
$(c_4)$ $1\to x = x = 1 \rightsquigarrow x.$\\
$(c_5)$ If $x\leq y$, then $z\to x\leq z\to y$ and
   $z\rightsquigarrow x\leq z\rightsquigarrow y.$\\
$(c_6)$ If $x\leq y$, then $y\to z\leq x\to z$ and
$y\rightsquigarrow  z\leq x\rightsquigarrow  z.$

\end{prop}

For all $x,y \in A$, define:
$$x \vee_1 y = (x \rightarrow y) \rightsquigarrow y, \:\:\: x \vee_2 y = (x \rightsquigarrow y) \rightarrow y.$$

\begin{prop} \label{psBCK-1.20-10} In any bounded pseudo-BCK algebra $A$ the following hold for all
$x,y \in A$:\\
$(1)$ $1 \vee_1 x = x \vee_1 1 = 1 = 1 \vee_2 x = x \vee_1 1.$ \\
$(2)$ $x \leq y$ implies $x \vee_1 y = y$ and $x \vee_2 y = y.$\\
$(3)$ $x \vee_1 x = x \vee_2 x = x.$\\
$(4)$ If $x_1 \le x_2$ and $y_1 \le y_2,$ then $x_1\vee_1 y_1 \le
x_2\vee_1 y_2$ and $x_1\vee_1 y_2 \le x_2\vee_1 y_2.$

\end{prop}
\begin{proof}
$(1)$ We have:
$1 \vee_1 x = (1 \rightarrow x) \rightsquigarrow x = 1$ and $x \vee_1 1 = (x \rightarrow 1) \rightsquigarrow 1 = 1$, so\\
$1 \vee_1 x = x \vee_1 1 = 1$, by $(c_4)$ of Proposition
\ref{psBCK-1.20-5}. Similarly, $1 \vee_2 x = x \vee_2 1 = 1$. \\
$(2)$ $x \vee_1 y = (x \rightarrow y) \rightsquigarrow y = 1 \rightsquigarrow y = y$. Similarly, $x \vee_2 y = y$.\\
$(3)$ By definitions of $\vee_1$ and $\vee_2$.\\
$(4)$ It is evident.
\end{proof}

\begin{prop} \label{psBCK-1.20-15}  In any bounded pseudo-BCK algebra $A$ the following hold for all
$x,y \in A$:\\
$(1)$ $x\vee_1 y^{- \sim}=x^{- \sim}\vee_1 y^{- \sim}$ and
      $x\vee_2 y^{\sim -}=x^{\sim -}\vee_2 y^{\sim -}.$ \\
$(2)$ $x\vee_1 y^{\sim}=x^{- \sim}\vee_1 y^{\sim}$ and
      $x\vee_2 y^-=x^{\sim -}\vee_2 y^-.$ \\
$(3)$ ${(x^{- \sim}\vee_1 y^{- \sim})}^{- \sim}= x^{- \sim}\vee_1 y^{- \sim}$ and
      ${(x^{\sim -}\vee_2 y^{\sim -})}^{\sim -}= x^{\sim -}\vee_2 y^{\sim -}$.
\end{prop}

\begin{proof}
The proof follows by direct computations.
\end{proof}

\begin{prop} \label{psBCK-1.20-20}
In any pseudo-BCK algebra the following hold for all $x,y \in A$:\\
$$x \vee_1 y \rightarrow y = x \rightarrow y\ \mbox{and}\ x \vee_2 y
\rightsquigarrow y = x \rightsquigarrow y. \eqno(2.1)$$
\end{prop}
\begin{proof} It is a consequence of property $(c_3)$.
\end{proof}

\begin{lemma} \label{psBCK-1.20-30} In any pseudo-BCK algebra $A$ we have\\
$(1)$ $x \vee_1 y \:\: (y \vee_1 x)$ is an upper bound of $\{x,y\}$.\\
$(2)$ $x \vee_2 y \:\: (y \vee_2 x)$ is an upper bound of $\{x,y\}.$
\end{lemma}
\begin{proof}
$(1)$ By $(A_2)$ we have $x \leq (x \rightarrow y) \rightsquigarrow y$. \\
Since by $(c_{2})$, $y \leq (x \rightarrow y) \rightsquigarrow y$, we conclude that $x,y \leq x \vee_1 y$.\\
Similarly we get $x,y \leq y \vee_1 x$.\\
$(2)$ Similarly as $(1)$.
\end{proof}

\begin{ex} \label{psBCK-1.25} (\cite{Ciu12})
Consider $A=\{o_1,a_1,b_1,c_1,o_2,a_2,b_2,c_2,1\}$ with
$o_1<a_1,b_1<c_1<1$ and $a_1, b_1$ incomparable, $o_2<a_2,b_2<c_2<1$
and $a_2, b_2$ incomparable. Assume also that any element of the set
$\{o_1,a_1,b_1,c_1\}$ is incomparable with any element of the set
$\{o_2,a_2,b_2,c_2\}$. Consider the operations
$\rightarrow,\rightsquigarrow$ given by the following tables:
\[
\hspace{5mm}
\begin{array}{c|ccccccccc}
\rightarrow
    & o_1 & a_1 & b_1 & c_1 & o_2 & a_2 & b_2 & c_2 & 1 \\ \hline
o_1 & 1   & 1   & 1   & 1   & o_2 & a_2 & b_2 & c_2 & 1 \\
a_1 & o_1 & 1   & b_1 & 1   & o_2 & a_2 & b_2 & c_2 & 1 \\
b_1 & a_1 & a_1 & 1   & 1   & o_2 & a_2 & b_2 & c_2 & 1 \\
c_1 & o_1 & a_1 & b_1 & 1   & o_2 & a_2 & b_2 & c_2 & 1 \\
o_2 & o_1 & a_1 & b_1 & c_1 & 1   & 1   & 1   & 1   & 1 \\
a_2 & o_1 & a_1 & b_1 & c_1 & o_2 & 1   & b_2 & 1   & 1 \\
b_2 & o_1 & a_1 & b_1 & c_1 & c_2 & c_2 & 1   & 1   & 1 \\
c_2 & o_1 & a_1 & b_1 & c_1 & o_2 & c_2 & b_2 & 1   & 1 \\
1   & o_1 & a_1 & b_1 & c_1 & o_2 & a_2 & b_2 & c_2 & 1
\end{array}
\hspace{10mm}
\begin{array}{c|ccccccccc}
\rightsquigarrow
    & o_1 & a_1 & b_1 & c_1 & o_2 & a_2 & b_2 & c_2 & 1 \\ \hline
o_1 & 1   & 1   & 1   & 1   & o_2 & a_2 & b_2 & c_2 & 1 \\
a_1 & b_1 & 1   & b_1 & 1   & o_2 & a_2 & b_2 & c_2 & 1 \\
b_1 & o_1 & a_1 & 1   & 1   & o_2 & a_2 & b_2 & c_2 & 1 \\
c_1 & o_1 & a_1 & b_1 & 1   & o_2 & a_2 & b_2 & c_2 & 1 \\
o_2 & o_1 & a_1 & b_1 & c_1 & 1   & 1   & 1   & 1   & 1 \\
a_2 & o_1 & a_1 & b_1 & c_1 & b_2 & 1   & b_2 & 1   & 1 \\
b_2 & o_1 & a_1 & b_1 & c_1 & b_2 & c_2 & 1   & 1   & 1 \\
c_2 & o_1 & a_1 & b_1 & c_1 & b_2 & c_2 & b_2 & 1   & 1 \\
1   & o_1 & a_1 & b_1 & c_1 & o_2 & a_2 & b_2 & c_2 & 1\end{array}
.
\]
Then $(A,\leq,\rightarrow,\rightsquigarrow,1)$ is a pseudo-BCK
algebra which is not a BCK algebra.
\end{ex}

\begin{ex} \label{psBCK-1.40} (\cite{Ciu12})
Consider $A=\{0,a,b,c,1\}$ with $0<a,b<c<1$ and $a, b$ incomparable.
Let the operations $\rightarrow,\rightsquigarrow$ be given by the
following tables:
\[
\hspace{10mm}
\begin{array}{c|ccccc}
\rightarrow & 0 & a & b & c & 1 \\ \hline
0 & 1 & 1 & 1 & 1 & 1 \\
a & 0 & 1 & b & 1 & 1 \\
b & a & a & 1 & 1 & 1 \\
c & 0 & a & b & 1 & 1 \\
1 & 0 & a & b & c & 1
\end{array}
\hspace{10mm}
\begin{array}{c|ccccc}
\rightsquigarrow & 0 & a & b & c & 1 \\ \hline
0 & 1 & 1 & 1 & 1 & 1 \\
a & b & 1 & b & 1 & 1 \\
b & 0 & a & 1 & 1 & 1 \\
c & 0 & a & b & 1 & 1 \\
1 & 0 & a & b & c & 1
\end{array}
.
\]
Then $(A,\leq,\rightarrow,\rightsquigarrow,0,1)$ is a bounded
pseudo-BCK algebra.
\end{ex}

\begin{Def} \label{psBCK-1.60} $\rm($\cite{Ior1}$\rm)$
A pseudo-BCK algebra with the \emph{$ (pP)$ condition} (i.e. with
the \emph{pseudo-product} condition) or  a \emph{pseudo-BCK$(pP)$
algebra} for short, is a pseudo-BCK algebra $A$ satisfying the $\rm
(pP)$ condition:\\  for all   $x,y \in A$, there exists $x\odot
y=\min\{z \mid x \leq y \rightarrow z\}=\min\{z \mid y \leq x
\rightsquigarrow z\}$.
\end{Def}

\begin{Def} \label{psBCK-1.63} $\rm($\cite{Ior1}$\rm)$
$(1)$ Let $(A,\leq,\rightarrow,\rightsquigarrow,1)$ be a pseudo-BCK
algebra. If the poset $(A,\leq)$ is a lattice, then we say that
$A$ is a \emph{pseudo-BCK lattice}.\\
$(2)$ Let $(A,\leq,\rightarrow,\rightsquigarrow,1)$ be a
pseudo-BCK$\rm (pP)$ algebra. If the poset $(A,\leq)$ is a lattice,
then we say that $A$ is a \emph{pseudo-BCK(pP) lattice}.
\end{Def}

\begin{Def} \label{psBCK-1.170} Let $A$ be a bounded pseudo-BCK algebra. Then:\\
$(1)$ $A$ is called \emph{good} if $x^{- \sim}=x^{\sim -}$ for all $x \in A$.\\
$(2)$ $A$ is with the \emph{$\rm (pDN)$ condition} (i.e. with the \emph{pseudo-double negation} condition) or a pseudo-BCK$\rm (pDN)$ algebra for short if
$x^{- \sim}=x^{\sim -}=x$ for all $x \in A$.
\end{Def}

For example, Examples \ref{psBCK-1.10-30} and \ref{psBCK-1.10-40}
are good pseudo-BCK algebras.

\begin{ex} \label{psBCK-1.190}
Consider the pseudo-BCK lattice $A$ from Example \ref{psBCK-1.40}.
Since $a^{- \sim}=1$ and $a^{\sim -}=b$, it follows that $A$ is not
good. $A$ can be embedded into the good pseudo-BCK lattice
$(A_1,\leq,\rightarrow,\rightsquigarrow,0,1)$, where
$A_1=\{0,a,b,c,d,1\}$ with $0 < a < b, c < d < 1$ and $b, c$ are
incomparable. The operations $\rightarrow$ and $\rightsquigarrow$
are defined as follows:
\[
\begin{array}{c|cccccc}
\rightarrow & 0 & a & b & c & d & 1 \\ \hline
0 & 1 & 1 & 1 & 1 & 1 & 1 \\
a & 0 & 1 & 1 & 1 & 1 & 1 \\
b & 0 & a & 1 & c & 1 & 1 \\
c & 0 & b & b & 1 & 1 & 1 \\
d & 0 & a & b & c & 1 & 1 \\
1 & 0 & a & b & c & d & 1
\end{array}
\hspace{10mm}
\begin{array}{c|cccccc}
\rightsquigarrow  & 0 & a & b & c & d & 1 \\ \hline
0 & 1 & 1 & 1 & 1 & 1 & 1 \\
a & 0 & 1 & 1 & 1 & 1 & 1 \\
b & 0 & c & 1 & c & 1 & 1 \\
c & 0 & a & b & 1 & 1 & 1 \\
d & 0 & a & b & c & 1 & 1 \\
1 & 0 & a & b & c & d & 1
\end{array}.
\]

One can easily check that $A_1$ is a good pseudo-BCK algebra.
Moreover, we can see that:\\
\hspace*{2cm} $\min\{z \mid c \leq b \rightarrow z\}=\min\{b,c,d,1\},$ \\
\hspace*{2cm} $\min\{z \mid b \leq c \rightsquigarrow z\}=\min\{b,c,d,1\}$ \\
do not exist. Thus, $c \odot b$ does not exist, so $ A_1$ is without the $\rm (pP)$ condition. \\
Since $(A_1,\leq)$ is a lattice, it follows that $ A_1$ is a good
pseudo-BCK lattice without the $\rm (pP)$ condition.
\end{ex}

\begin{Def} \label{psBCK-1.220} $\rm($\cite{Geo15}, \cite{Ior13}$\rm)$ Let $A$ be a pseudo-BCK algebra. \\
$(1)$ If $x \vee_1 y = y \vee_1 x$ for all $x,y \in A$, then $A$ is called \emph{$\vee_1$-commutative}. \\
$(2)$ If $x \vee_2 y = y \vee_2 x$ for all $x,y \in A$, then $A$ is
called \emph{$\vee_2$-commutative}.
\end{Def}

\begin{lemma} \label{psBCK-1.230} $\rm($\cite{Geo15}, \cite{Ior13}$\rm)$ If $A$ is a pseudo-BCK algebra, then:\\
$(1)$ $A$ is $\vee_1$ -commutative if and only if it is a
join-semilattice (under $\leq$).\\
$(2)$ $A$ is $\vee_2$ -commutative if and only  if it is a
join-semilattice (under $\leq$).
\end{lemma}

\begin{Def} \label{psBCK-1.240} $\rm($\cite{Geo15}, \cite{Ior13}$\rm)$
A pseudo-BCK algebra is called \emph{sup-commutative} if it is both
$\vee_1$-commutative and $\vee_2$-commutative.
\end{Def}

\begin{theo} \label{psBCK-1.250} $\rm($\cite{Geo15}, \cite{Ior13}$\rm)$
A pseudo-BCK algebra is sup-commutative if and only if it is a
semilattice with respect to both $\vee_1$ and $\vee_2$.
\end{theo}

\begin{cor} \label{psBCK-1.260} $\rm($\cite{Geo15}, \cite{Ior13}$\rm)$
If $A$ is a sup-commutative pseudo-BCK algebra, then $\vee_1=\vee_2$
for all $x,y \in A$, so $A$ is a \emph{join-semilattice} with
$x\vee_1 y=x \vee_2 y$.
\end{cor}

If $A$ is a sup-commutative pseudo-BCK algebra, then $\vee_1=\vee_2$
for all $x,y \in A$, so $A$ is a semilattice with $x\vee y=x\vee_1
y=x \vee_2 y$.
\begin{proof}
It follows applying \cite[Cor.1.17]{Geo15}.
\end{proof}

In what follows, we introduce a new class of pseudo-BCK algebras
which will be used in the next sections.

\begin{Def} \label{psBCK-1.270}
A bounded pseudo-BCK algebra $A$ is said to be with the
\emph{join-negation} (JN for short) if
$$x\vee_1 y= x^{- \sim}\vee_1 y^{- \sim}\quad  \mbox{and}\quad
                x\vee_2 y= x^{- \sim}\vee_2 y^{- \sim}\quad  \mbox{for all}\
                x,y\in A.$$
\end{Def}

\begin{rem} \label{psBCK-1.280}
$(1)$ Every bounded pseudo-BCK(pDN) algebra is with (JN). \\
$(2)$ Every bounded sup-commutative pseudo-BCK algebra is with (pDN),
so it is with (JN). \\
$(3)$ Every locally finite pseudo-hoop is a bounded pseudo-BCK(pDN) algebra,
so it is with (JN) (see \cite{Ciu14}).
\end{rem}

We recall that a \emph{downwards-directed set} (or a \emph{filtered
set}) is a partially ordered set $(A,\leq)$ such that whenever
$a,b\in A,$
there exists $x\in A$ such that $x\leq a$ and $x\leq b$.\\

According to \cite{DvVe}, we say that a pseudo-BCK algebra $A$
satisfies the \emph{relative cancellation property}, (RCP) for
short, if for every $a,b,c\in A$,
$$
a,b\leq c\quad \mbox{and}\quad  c\rightarrow a=c\rightarrow b,\ \
c\squig a=c\squig b  \quad \mbox{imply}\quad a=b.
$$

We note that a pseudo-BCK algebra $A$ that is sup-commutative and
satisfies the (RCP)-condition is said to be a {\it \L ukasiewicz
pseudo-BCK algebra,} see \cite{DvVe}.

\begin{ex} \label{psBCK-1.290}
The pseudo-BCK algebra $A$ from Example \ref{psBCK-1.40} is
downwards-directed with (RCP).
\end{ex}

\begin{prop} \label{psBCK-1.300}
Any downwards-directed sup-commutative pseudo-BCK algebra has (RCP).
\end{prop}
\begin{proof}
Consider $a,b,c\in A$ such that $a,b\leq c$. There exists $x\in A$ such that $x\leq a,b$.\\
By $(c_6)$, from $a\leq c$ it follows that $c\rightsquigarrow x\leq a\rightsquigarrow x$.\\
According to Proposition \ref{psBCK-1.20-10}$(2)$ and $(c_1)$ we have:
\begin{eqnarray*}
a\rightsquigarrow x&=&(c\rightsquigarrow x)\vee_1(a\rightsquigarrow
x)=(a\rightsquigarrow x)\vee_1(c\rightsquigarrow x)\\
&=&
  [(a\rightsquigarrow x)\rightarrow (c\rightsquigarrow x)]\rightsquigarrow
  (c\rightsquigarrow x)=
  c\rightsquigarrow [(a\rightsquigarrow x)\rightarrow x]
  \rightsquigarrow (c\rightsquigarrow x)\\&=&
  [c\rightsquigarrow (a\vee_2 x)]\rightsquigarrow (c\rightsquigarrow x)=
  (c\rightsquigarrow a)\rightsquigarrow (c\rightsquigarrow x).
\end{eqnarray*}
Similarly, $b\rightsquigarrow x = (c\rightsquigarrow b)\rightsquigarrow (c\rightsquigarrow x)=a\rightsquigarrow x$.\\
We have: $\:$
$a=x\vee_2a=a\vee_2 x=(a\rightsquigarrow x)\rightarrow x=(b\rightsquigarrow x)
\rightarrow x=b\vee_2 x=x\vee_2b=b$. \\
Thus, $A$ has (RCP).
\end{proof}

We say that a nonempty subset $F$ of a pseudo-BCK algebra $A$ is a
{\it filter} (or a {\it deductive system}, \cite{Kuhr, Kuhr1}) if
(i) $1 \in F$, and (ii) if $a \in F$ and $a\to b \in F,$ then $b \in
F.$ It is easy to verify that a set $F$ containing $1$ is a filter
if and only if (ii)' if $a \in F$ and $a\rightsquigarrow b \in F,$
then $b \in F.$  A filter $F$ is  (i) {\it maximal} if it is a
proper subset of $A$ and not properly contained in another proper
filter of $A,$ (ii) {\it normal} if $a\to b \in F$ if and only if
$a\rightsquigarrow b\in F.$  Given a normal filter $F$, the relation
 $\Theta_F$ on $A$ given by
$$
(a,b) \in \Theta_F\ \Leftrightarrow a\to b\in F\ \mbox{and}\ b\to
a\in F
$$
is a congruence. Then $F = [1]\Theta_F$ and   the quotient class
$A/F$ defined as $A/\Theta_F$  is again a pseudo-BCK algebra, and we
write $a/F=[a]_{\Theta_F}$ for every $a\in A,$ see \cite{Kuhr,
Kuhr1}.

Given an integer $n\ge 1,$ we define inductively

$$x\to^0y =y,\quad x\to^n y = x\to (x\to^{n-1} y), \quad n\ge 1,$$
and
$$x\rightsquigarrow^0y =y,\quad x\rightsquigarrow^n y = x\rightsquigarrow
 (x\rightsquigarrow^{n-1} y), \quad n\ge 1.$$

For any nonempty system $X$ of a pseudo-BCK algebra $A,$ there is
the least filter of $A$ generated by $X,$ we denote it by $F(X);$ in
particular, if $X=\{u\}$ is a singleton, we set $F(u):=F(\{u\}).$ By
\cite{Kuhr, Kuhr1},  we have

$$F(X)=\{a\in A:\ x_1\to(\cdots \to (x_n\to a)\cdots )=1,\
x_1,\ldots,x_n \in X,\ n\ge 1\}\eqno(2.2)
$$
and
$$F(u)=\{a\in A:\ u\to^n a=1,\ n \ge 1\}= \{a\in A:\ u\rightsquigarrow^n a=1,\ n \ge 1\}.\eqno(2.3)$$

If $F$ is a filter and $b\in A$, then the filter, $F_b,$ of $A$
generated by $F\cup\{b\}$ is the set
$$ F_b=\{a\in A:\ b\to ^n a\in F\ \mbox{for some}\ n \in \mathbb N\}.\eqno(2.4)$$

\section{States on Pseudo-BCK Algebras}

We present a notion of states on bounded pseudo-BCK algebras. We
characterize extremal states as state-morphisms and we show that the
quotient pseudo-BCK algebra through the kernel of a state is always
an MV-algebra. We emphasize  that our characterizations of states
can be studied without the assumption that the pseudo-BCK algebra is
a pseudo MV-algebra or a $\vee$-lattice, see \cite{Dvu2, Kuhr}.

\begin{Def}\label{psBCK-2.10}
A {\it Bosbach state} on a bounded pseudo-BCK algebra $A$ is a function $s:A\longrightarrow [0,1]$ such that the following conditions hold for any $x,y\in A:$\\
$(B_1)$ $s(x)+s(x\rightarrow y)=s(y)+s(y\rightarrow x);$\\
$(B_2)$ $s(x)+s(x\rightsquigarrow y)=s(y)+s(y\rightsquigarrow x);$\\
$(B_3)$ $s(0)=0$ and $s(1)=1.$
\end{Def}

\begin{ex} \label{psBCK-2.20}
Consider the bounded pseudo-BCK lattice $A$ from Example
\ref{psBCK-1.190}. The function $s : A \longrightarrow [0,1]$
defined by: $s(0)=0, s(a)=1, s(b)=1, s(c)=1, s(d) = 1, s(1)=1$ is a
unique Bosbach state on $A$.
\end{ex}

Not every bounded pseudo-BCK algebra has a Bosbach state:

\begin{ex} \label{psBCK-2.25}
Consider the bounded pseudo-BCK lattice  $A$ from Example
\ref{psBCK-1.40}. One can prove that $A$ has no Bosbach state.
Indeed, assume that $A$ admits a Bosbach state $s$ such that
$s(0)=0$, $s(a)=\alpha$, $s(b)=\beta$, $s(c)=\gamma$, $s(1)=1$.
From $s(x)+s(x\rightarrow y)=s(y)+s(y\rightarrow x)$, taking $x=a, y=0$, $x=b,y=0$ and respectively $x=c,y=0$ we get $\alpha=1$, $\beta = 0$, $\gamma = 1$.\\
On the other hand, taking $x=b, y=0$ in $s(x)+s(x\rightsquigarrow
y)=s(y)+s(y\rightsquigarrow x)$ we get $\beta + 0 = 0 + 1$, so $0=1$
which is a contradiction. Hence, $A$ does not admit a Bosbach state.
\end{ex}

\begin{prop}\label{psBCK-2.275-10}
Let $A$ be a bounded pseudo-BCK algebra and $s$ a Bosbach state on
$A$. For all $x,y\in A,$ the following properties hold:\\
$(1)$ $s(y\to x)=1+s(x)-s(y) =s(y\rightsquigarrow x)$ and $s(x)\le
s(y)$ whenever $x\le
y.$\\
$(2)$ $s(x\vee_1 y) = s(y \vee_1 x)$ $\:$ and $\:$ $s(x\vee_2 y) = s(y \vee_2 x).$\\
$(3)$ $s(x\vee_1 y^{- \sim})=s(x^{- \sim}\vee_1 y^{- \sim})$ and
      $s(x\vee_2 {y^{\sim}}^-)=s({x^{\sim}}^-\vee_2 y^{\sim -}).$  \\
$(4)$ $s(x^{- \sim}\vee_1 y)=s(x\vee_1 y^{- \sim})$ and
      $s(x^{\sim -}\vee_2 y)=s(x\vee_2 y^{\sim -}).$\\
$(5)$ $s(x^{-\sim}) = s(x) = s(x^{\sim-}).$\\
$(6)$ $s(x^-)= 1-s(x)=s(x^\sim).$
\end{prop}

\begin{proof} (1) It is straightforward.

(2) By (2.1) and property (1), we have $ s(x\to y)= s(x\vee_1 y \to
y) = 1 +s(y)-s(x\vee_1 y)$ and $s(y\to x)=s(y\vee_1 x\to x) = 1
+s(x)-s(y\vee_1 x).$ Then $s(x\to y) = s(y) +s(y\to x) -s(x)=
s(y)+(1+s(x) -s(y\vee_1x))-s(x)$ proving $s(x\vee_1 y)=s(y\vee_1
x).$ Similarly, $s(x\vee_2 y)=s(y\vee_2 x).$

(3) and (4) follow from Proposition \ref{psBCK-1.20-15}.

(5) Since  $x^{-\sim} = x\vee_1 0,$ by (2) we have $s(x^{-\sim}) =
s(x\vee_1 0) = s(0\vee_1 x) = s((0\to x)\squig x) = s(x).$  In a
similar way, we have $s(x) = s(x^{\sim-}).$

(6) $s(x^-)=s(x\to 0) = s(0)-s(x)+s(0\to x) = 1-s(x).$
\end{proof}

\begin{prop}\label{psBCK-2.275-20} Let $A$ be a bounded pseudo-BCK algebra
and $s$ a Bosbach state on $A$. For all $x,y\in A,$ the following properties hold:\\
$(1)$ $s(x^{- \sim}\rightarrow y)=s(x\rightarrow y^{- \sim})$ and
      $s(x^{\sim -}\rightsquigarrow y)=s(x\rightsquigarrow y^{\sim -}).$ \\
$(2)$ $s(x\rightarrow y^{- \sim})=s(y^-\rightsquigarrow x^-)=s(x^{-
\sim}\rightarrow y^{- \sim})=
       s(x^{- \sim}\rightarrow y)$ and \\
$\hspace*{0.5cm}$
 $s(x\rightsquigarrow y^{\sim -})=s(y^{\sim}\rightarrow x^{\sim})=s(x^{\sim -}\rightsquigarrow y^{\sim -})=
  s(x^{\sim -}\rightsquigarrow y).$ \\
$(3)$ $s(x^{\sim}\rightarrow y^{- \sim})=s(x^{\sim}\rightarrow y)$
and
      $s(x^-\rightsquigarrow y^{\sim -})=s(x^-\rightsquigarrow y)$.
\end{prop}
\begin{proof}
$(1)$ Using Proposition \ref{psBCK-2.275-10}(4), we have:\\
$s(x^{- \sim}\rightarrow y)=1-s(x^{- \sim}\vee_1 y)+s(y)=1-s(x\vee_1
y^{- \sim})+s(y^{- \sim})=
s(x\rightarrow y^{- \sim})$.\\
$(2)$ It follows by $(c_4)$ and $(1)$. \\
$(3)$ Applying Proposition \ref{psBCK-2.275-10}$(4)$ we get:\\
$s(x^{\sim}\rightarrow y^{- \sim})=1-s(x^{\sim}\vee_1 y^{-
\sim})+s(y^{- \sim})=
1-s(x^{\sim - \sim}\vee_1 y)+s(y)=\\
1-s(x^{\sim}\vee_1 y)+s(y)=s(x^{\sim}\rightarrow y)$.\\
Similarly, $s(x^-\rightsquigarrow
{y^{\sim}}^-)=s(x^-\rightsquigarrow y)$.
\end{proof}

\begin{prop} \label{psBCK-2.280} Let $A$ be a bounded pseudo-BCK algebra and a function
$s : A\longrightarrow [0,1]$ such that $s(0)=0$, $s(x\vee_1 y)=s(y\vee_1 x)$ and
$s(x\vee_2 y)=s(y\vee_2 x)$ for all $x,y \in A$. Then the following are equivalent:\\
$(a)$ $s$ is a Bosbach state on $A.$\\
$(b)$ for all $x,y \in A$, $y \leq x$ implies $s(x \rightarrow y)= s(x \rightsquigarrow y)=1-s(x)+s(y).$\\
$(c)$ for all $x,y \in A$,
$s(x \rightarrow y)= 1-s(x \vee_1 y) + s(y)$ and \\
$\hspace*{3.5cm}$ $s(x \rightsquigarrow y)= 1-s(x \vee_2 y) + s(y)$.
\end{prop}
\begin{proof}
$(a)\Rightarrow (b).$ It follows from Proposition
\ref{psBCK-2.275-10}(1).\\
$(b)\Rightarrow (c).$ It follows from the proof of Proposition
\ref{psBCK-2.275-10}(2).\\
$(c)\Rightarrow (a).$ Using Proposition \ref{psBCK-2.275-10}(2), we
get:
\begin{eqnarray*}
s(x)+s(x \rightarrow y)&=&s(x)+1-s(x \vee_1 y)+s(y)\\
&=&1-s(y \vee_1 x)+s(x)+s(y)=s(y)+s(y \rightarrow x).
\end{eqnarray*}
Similarly,
\begin{eqnarray*}
s(x)+s(x \rightsquigarrow y)&=&s(x)+1-s(x \vee_2 y)+s(y)\\
&=&1-s(y \vee_2 x)+s(x)+s(y)=s(y)+s(y \rightsquigarrow x).
\end{eqnarray*}
Moreover, by $(c)$ we have:\\
$\hspace*{3cm}$
$s(1)=s(x \rightarrow x)=1-s(x)+s(x)=1$.\\
Thus, $s$ is a Bosbach state on $A$.
\end{proof}

The following proposition is crucial for our study.

\begin{prop}\label{AD4}  Let $s$ be a Bosbach state on a bounded pseudo-BCK
algebra.  Then, for all $x,y \in A,$ we have:\\
$(1)$ $s(x\vee_1 y) = s(x\vee_2 y).$\\
$(2)$ $s(x\to y) = s(x\squig y)$.
\end{prop}

\begin{proof} (1) First we prove the equality for $y \le x.$

Using Proposition \ref{psBCK-2.275-10}(2), we have $s(x\vee_1 y)=
s(y\vee_1 x) = s(x)$ and by Proposition \ref{psBCK-1.20-10}(2)
$s(x\vee_2 y) = s(x),$ i.e., $s(x\vee_1 y)=s(x\vee_2y).$

Assume now that $x$ and $y$ are arbitrary elements of $A$. Using
again  Proposition \ref{psBCK-2.275-10}(2) and the first part of the
proof, we have
\begin{eqnarray*}
s(x\vee_1y)&=& s(x\vee_1 (x\vee_1y)) = s((x\vee_1 y)\vee_1 x)\\
&=& s((x\vee_1 y)\vee_2 x) \ge s(y\vee_2 x)\\
&=& s(x \vee_2 (y\vee_2 x) ) \ge s(x\vee_2 y)\\
&=&    s(y\vee_2 (x \vee_2y))=  s((x\vee_2 y) \vee_2 y)\\
&\ge& s(x\vee_1 y).
\end{eqnarray*}

(2) This follows immediately from Proposition \ref{psBCK-2.280}(c)
and the first equation.
\end{proof}

Consider the real interval $[0,1]$ of reals equipped with the \L
ukasiewicz implication
$\rightarrow_{\mbox{\tiny \L}}$ defined by\\
$$x \rightarrow_{\mbox{\tiny \L}} y = x^-\oplus y=\min\{1-x+y,1\},\quad x,y
\in [0,1].$$

\begin{Def}\label{psBCK-2.280-10} Let $A$ be a bounded
pseudo-BCK algebra. A \emph{state-morphism} on
$A$ is a function $m : A\longrightarrow [0,1]$ such that:\\
$(SM_1)$ $m(0)=0.$\\
$(SM_2)$ $m(x \rightarrow y) = m(x \rightsquigarrow y)=m(x)
\rightarrow_{\mbox{\tiny \L}} m(y),$ $x,y \in A.$
\end{Def}

\begin{prop}\label{psBCK-2.280-20} Every state-morphism on a bounded
pseudo-BCK algebra $A$ is a Bosbach state on $A$.
\end{prop}
\begin{proof} It is obvious that $m(1)=m(x \rightarrow x)=m(x) \rightarrow_{\mbox{\tiny \L}} m(x)=1$. \\
We also have:
\begin{eqnarray*}
m(x)+m(x \rightarrow y)&=&m(x)+m(x) \rightarrow_{\mbox{\tiny \L}} m(y) =m(x)+\min\{1-m(x)+m(y),1\} \\
&=&\min\{1+m(y),1+m(x)\}=m(y)+\min\{1-m(y)+m(x),1\} \\
&=&m(y)+m(y) \rightarrow_{\mbox{\tiny \L}} m(x)=m(y)+m(y \rightarrow
x).
\end{eqnarray*}
Similarly, $m(x)+m(x \rightsquigarrow y)= m(y) + m(y \rightsquigarrow x)$.\\
Thus, $s$ is a Bosbach state on $A$.
\end{proof}

\begin{prop}\label{psBCK-2.280-30} Let $A$ be a
bounded pseudo-BCK algebra. A Bosbach state $m$ on $A$ is a
state-morphism if and only if:
$$m(x \vee_1 y)=\max\{m(x),m(y)\}$$
for all $x,y \in A,$ or equivalently,
$$m(x \vee_2
y)=\max\{m(x),m(y)\}$$ for all $x,y \in A.$
\end{prop}
\begin{proof}  In view of Proposition \ref{AD4}, the two equations
are equivalent. If $m$ is a state-morphism on $A$, then by
Proposition \ref{psBCK-2.280-20} $m$ is a Bosbach state.
Using the relation:\\
$m(x \rightarrow y)=1-m(x \vee_1 y)+m(y),$ we obtain:
\begin{eqnarray*}
m(x \vee_1 y) &=& 1 + m(y) - m(x \rightarrow y)=1+m(y)-(m(x)
\rightarrow_{\mbox{\tiny \L}} m(y))\\
 &=& 1+m(y)-\min\{1-m(x)+m(y),1\}\\
 &=&
1+m(y)+\max\{-1+m(x)-m(y),-1\}=\max\{m(x),m(y)\}.
\end{eqnarray*}
For the converse, assume that $m$ is a Bosbach state on $A$ such that \\
$\hspace*{2cm}$ $m(x \vee_1 y)=\max\{m(x),m(y)\}$ for all $x,y \in A$.\\
Then, using again the relation:\\
$m(x \rightarrow y)=1-m(x \vee_1 y)+m(y),$ we have :
\begin{eqnarray*}
m(x \rightarrow y) &=& 1+m(y)-\max\{(m(x),m(y)\}\\
&=&1+m(y)+\min\{-m(x),-m(y)\}\\
&=&\min\{1-m(x)+m(y),1\}= m(x)\rightarrow_{\mbox{\tiny \L}} m(y).
\end{eqnarray*}
Similarly,\\
$m(x \rightsquigarrow y) = m(x)\rightarrow_{\mbox{\tiny \L}} m(y)$.\\
Thus, $m$ is a state-morphism on $A$.
\end{proof}

Now we present an example of a  linearly ordered pseudo-BCK algebra
without the (pP) condition but having a unique Bosbach state (= a
unique state-morphism). On the other hand, not every linearly
ordered pseudo-BCK algebra admits a Bosbach state, see Example
\ref{psBCK-2.390}.

\begin{ex} \label{psBCK-2.280-40} (\cite{Ciu12})
Consider $A=\{0,a,b,c,1\}$ with $0<a<b, c<1$ and $b, c$
incomparable. Consider the operations $\rightarrow,\rightsquigarrow$
given by the following tables:
\[
\hspace{10mm}
\begin{array}{c|ccccc}
\rightarrow & 0 & a & b & c & 1 \\ \hline
0 & 1 & 1 & 1 & 1 & 1 \\
a & 0 & 1 & 1 & 1 & 1 \\
b & 0 & a & 1 & c & 1 \\
c & 0 & b & b & 1 & 1 \\
1 & 0 & a & b & c & 1
\end{array}
\hspace{10mm}
\begin{array}{c|ccccc}
\rightsquigarrow & 0 & a & b & c & 1 \\ \hline
0 & 1 & 1 & 1 & 1 & 1 \\
a & 0 & 1 & 1 & 1 & 1 \\
b & 0 & c & 1 & c & 1 \\
c & 0 & a & b & 1 & 1 \\
1 & 0 & a & b & c & 1
\end{array}
.
\]
Then $(A,\leq,\rightarrow,\rightsquigarrow,0,1)$ is a
bounded pseudo-BCK algebra.\\
Since $(A,\leq)$ is a lattice, it follows that $A$ is a pseudo-BCK lattice.\\
Moreover, we can see that $c \odot b = \min\{z \mid c \leq b
\rightarrow z\}=\min\{b,c,1\}$ does not exist. Hence, $A$ is a
pseudo-BCK lattice without the (pP) condition
 and the function
$m : A \longrightarrow [0,1]$ defined by: \\
$\hspace*{2.5cm}$ $m(0)=0, m(a)=1, m(b)=1, m(c)=1, m(1)=1$ \\
is a unique state-morphism on $A$.\\
Moreover, $m(x \vee_1 y)=m(x \vee_2 y)=\max\{m(x),m(y)\}$ for all
$x,y \in A$, hence $m$ is also a Bosbach state on $A$.
\end{ex}

The set $$\Ker(s):=\{a\in A\mid s(a)=1\}$$ is called the
\emph{kernel} of a Bosbach state $s$ on $A.$

\begin{prop} \label{psBCK-2.290}
Let $A$ be a bounded pseudo-BCK algebra and let $s$ be a Bosbach
state on $A$. Then $\Ker(s)$ is a proper and normal filter of $A$.
\end{prop}
\begin{proof} Obviously, $1\in \Ker(s)$ and $0\notin \Ker(s)$. \\
Assume that $a,a\rightarrow b \in \Ker(s)$. We have $1=s(a)\leq
s(a\vee_1 b)$, so $s(a\vee_1 b)=1$.\\
It follows that $1=s(a\rightarrow b)=1-s(a\vee_1 b)+s(b)=s(b)$.\\
Hence $b\in \Ker(s)$, so $\Ker(s)$ is a proper filter of $A$.\\
By Proposition \ref{AD4}(2),  $s(a\rightsquigarrow b)=s(a\rightarrow
b)$, and this proves that $\Ker(s)$ is normal.
\end{proof}

\begin{lemma} \label{psBCK-2.300}
Let $s$ be a Bosbach state on a bounded pseudo-BCK algebra $A$ and
$K=\Ker(s)$. In the bounded quotient pseudo-BCK  algebra
$(A/K,\leq, \rightarrow,\rightsquigarrow,0/K,1/K)$ we have:\\
$(1)$ $a/K \leq b/K$ iff $s(a\rightarrow b)=1$ iff $s(a\vee_1 b)=s(b)$ iff $s(a\vee_2 b)=s(b).$ \\
$(2)$ $a/K=b/K$ iff $s(a\rightarrow b)=s(b\rightarrow a)=1$ iff
$s(a)=s(b)=s(a\vee_1 b)$ iff $s(a\rightsquigarrow
b)=s(b\rightsquigarrow a)=1$ iff $s(a)=s(b)=s(a\vee_2 b)$.

Moreover, the mapping $\hat s:\ A/K \to [0,1]$ defined by $\hat
s(a/K):=s(a)$ $(a\in A)$  is a Bosbach state on $A/K.$
\end{lemma}

\begin{proof}
$(1)$ It follows easily: $a/K\leq b/K$ iff $(a\rightarrow
b)/K=a/K\rightarrow b/K=1/K=K$ iff
$a\rightarrow b\in K$ iff $s(a\rightarrow b)=1$.\\
As $s(a\rightarrow b)=1-s(a\vee_1 b)+s(b)$, we get $a/K\leq b/K$ iff $s(a\vee_1 b)=s(b).$\\
Similarly, $a/K\leq b/K$ iff $(a\rightsquigarrow
b)/K=a/K\rightsquigarrow b/K=1/K=K$ iff
$a\rightsquigarrow b\in K$ iff $s(a\rightsquigarrow b)=1$.\\
As $s(a\rightsquigarrow b)=1-s(a\vee_2 b)+s(b)$, we get $a/K\leq b/K$ iff $s(a\vee_2 b)=s(b).$\\
$(2)$ It follows easily from $(1)$.

The fact that $\hat s$ is a well-defined  Bosbach state on $A/K$ is
now straightforward.
\end{proof}

\begin{prop}\label{AD1} Let $s$ be a Bosbach state on a bounded
pseudo-BCK algebra $A$ and let $K=\Ker(s).$ For every element $x\in
A,$ we have
$$x^{-\sim}/K = x/K = x^{\sim-}/K,$$
that is, $A/K$ satisfies the (pDN) condition.
\end{prop}

\begin{proof} On one side, we have $x \le x^{-\sim}.$ On the other
one, by definition of a Bosbach state and Proposition
\ref{psBCK-2.275-10}(5), we have $s(x^{-\sim} \to x) = s(x) + s(x\to
x^{-\sim})-s(x^{-\sim}) = s(x\to x^{-\sim}) = s(1) = 1.$ Hence,
$x^{-\sim}/K = x/K.$  In a similar way, we prove the second
identity.
\end{proof}

\begin{rem}\label{AD3}  Let $s$ be a Bosbach state on a pseudo-BCK algebra $A$. According to
the proof of Proposition \ref{AD1}, we have $s(x^{-\sim}\to
x)=1=s(x^{\sim-}\to x)$ and $s(x^{-\sim}\squig
x)=1=s(x^{\sim-}\squig x).$
\end{rem}

\begin{prop}\label{AD2}  Let $s$ be a Bosbach state on a bounded
pseudo-BCK algebra $A.$ Then $A/K$ is $\vee_1$-commutative as well
as $\vee_2$-commutative, where $K=\Ker(s).$ In addition, $A/K$ is a
$\vee$-semilattice and good.
\end{prop}

\begin{proof}  Since $s$ is a Bosbach state, $A/K$ is a
BCK algebra. We denote by $\bar{x}:= x/K$, $x \in A$ and $\hat
s(a):=s(a)$ $(a \in A)$ is a Bosbach state on $A/K.$

(1) We show that if $\bar{x} \le \bar{y},$ then
$$ \bar{x}\vee_1 \bar{y} = \bar{y} = \bar{y} \vee_1 \bar{x}.\eqno(**)$$

By Proposition \ref{psBCK-1.20-10}, we have $\bar{x} \vee_1 \bar{y}
= \bar{y}.$  We have to show that $s((y\vee_1 x)\to y)=1.$
Calculate: By Proposition \ref{psBCK-2.275-10}(1), we have $s(y
\vee_1 x)= s((y\to x)\squig x)= \hat s((\bar y\to \bar x)\squig \bar
x)=  1 +\hat s(\bar x) - \hat s(\bar y\to \bar x)= 1 + \hat s(\bar
x) -[1+\hat s(\bar x)-\hat s(\bar y)] = \hat s(\bar y) = s(y).$\\
Therefore, using Proposition \ref{psBCK-2.275-10}(2), \\ $s((y\vee_1
x)\to y) = \hat s((\bar y \vee_1 \bar x)\to \bar y)= 1 +\hat s(\bar
y) - \hat s(\bar y\vee_1 \bar x) =  1 +\hat s(\bar y) - \hat s(\bar
x\vee_1 \bar y) =  1 +\hat s(\bar y) - \hat s(\bar y) = 1.$ Hence,
$(**)$ holds for $\bar x\le \bar y.$

(2) Now we show that $(**)$ holds for all $x,y \in A.$  By (1), we
have

$\bar x \vee_1 \bar y = \bar x \vee_1 (\bar x\vee_1 \bar y) = (\bar
x\vee_1 \bar y)\vee_1 \bar x \ge \bar y \vee_1 \bar x = \bar y
\vee_1 (\bar y\vee_1 \bar x)=(\bar y \vee_1 \bar x) \vee_1 \bar y
\ge \bar x \vee_1 \bar y.$

This implies that $A/K$ is $\vee_1$-commutative. In a similar way we
prove that $A/K$ is  $\vee_2$-commutative.  By Lemma
\ref{psBCK-1.20-30}, $A/K$ is a $\vee$-semilattice.
\end{proof}

\begin{prop} \label{psBCK-2.420} $\rm($\cite{Ciu12}$\rm)$
Let $A$ be a bounded good pseudo-BCK algebra. We define a binary
operation $\oplus$ on $A$ by $x \oplus y := x^{\sim} \rightarrow
{y^{\sim}}^-$. For all $x,y \in A,$ the
following hold:\\
$(1)$ $x \oplus y = y^- \rightsquigarrow x^{\sim -}$, \\
$(2)$ $x,y \leq x \oplus y$, \\
$(3)$ $x \oplus 0=0 \oplus x = x^{\sim -}$,\\
$(4)$ $x \oplus 1= 1 \oplus x = 1$,\\
$(5)$ $(x \oplus y)^{- \sim}=x \oplus y=x^{- \sim} \oplus y^{- \sim}$,\\
$(6)$ $\oplus$ is associative.
\end{prop}

An MV-algebra is an algebra $(A,\oplus,\odot,^-,0,1)$ of type
$\langle 2,2, 1,0,0\rangle$ such that (i) $\oplus$ is commutative
and associative, (ii) $0^-=1,$ (iii) $x\oplus 0 = x,$ (iv) $x \oplus
1 = 1,$ (v) $x^{**}= x,$ (vi) $y \oplus (y \oplus x^-)^- = x \oplus
(x\oplus y^-)^-$, and (vi) $x\odot y = (x^-\oplus y^-)^-.$  If we
define $x\to y =x\rightsquigarrow y = x^-\oplus y$, then
$(A,\to,\rightsquigarrow,0,1)$ is a bounded pseudo-BCK algebra.

An {\it MV-state} on an MV-algebra $A$ is a mapping $s:A\to [0,1]$
such that $s(1)=1$ and $s(a\oplus b) =s(a)+s(b)$ whenever $a\odot
b=0.$ Every MV-algebra admits at least one MV-state, and due to
\cite{Geo1}, every MV-state on $A$ coincides with a Bosbach state on
the BCK algebra $A$ and vice versa.

We note that the {\it radical}, $\Rad(A),$ of an MV-algebra $A$ is
the intersection of all maximal ideals of $A,$ \cite{CDM}.

\begin{prop}{\rm (\cite[Thm 6.1.32]{DvPu})} \label{arh-mv-30}
In any MV-algebra $A$ the following conditions are equivalent:\\
$(a)$ $\Rad(A)=\{0\}.$\\
$(b)$ $nx\leq x^{-}$ for all $n\in\NN$ implies $x=0.$\\
$(c)$ $nx\leq y^{-}$ for all $n\in \NN$ implies $x\wedge y=0.$\\
$(d)$ $nx\leq y$ for all $n\in\NN$ implies $x\odot y=x,$\\
where $nx=x_1\oplus \cdots \oplus x_n$ with $x_1=\cdots =x_n=x.$
\end{prop}

\begin{rem} \label{arh-mv-60}
An MV-algebra $A$ is archimedean in the  sense of \cite{DvPu} if it
satisfies the condition $(b)$ of Proposition \ref{arh-mv-30} and $A$
is archimedean in Belluce's sense \cite{Bell1} if it satisfies
the condition $(d)$ of Proposition \ref{arh-mv-30}.\\
By Proposition \ref{arh-mv-30} the two definitions of archimedean
MV-algebras are equivalent.
\end{rem}

\begin{theo} \label{psBCK-2.310} Let $s$ be a Bosbach state on a pseudo-BCK
algebra $A$ and let $K=\Ker(s)$. Then $(A/K,\oplus,^-,0/K)$, where
$$
a/K \oplus b/K=(a^{\sim} \rightarrow b)/K\ \mbox{and} \
(a/K)^-=a^-/K,
$$
is an archimedean MV-algebra and the map $\hat{s}(a/K):=s(a)$ is an
MV-state on this MV-algebra.
\end{theo}

\begin{proof}  By Propositions \ref{AD1} and \ref{AD2},
$A/K$ is a good pseudo-BCK algebra that is a $\vee$-semilattice and
$\hat s$ on $A/K$ is a Bosbach state such that $\Ker(\hat s)=
\{1/K\}.$ Due to \cite[Prop 3.4.7]{Kuhr}, $(A/K)/\Ker(\hat s)$ is
term-equivalent to an MV-algebra that is archimedean and $\hat s$ is
an MV-state on it. Since  $A/K = (A/K)/\Ker(\hat s),$ the same is
true also for $A/K,$ and this proves the theorem.
\end{proof}

We recall that if a pseudo-BCK algebra $A$ is good, in view of
Proposition \ref{psBCK-2.420}, we can define a binary operation
$\oplus$ via $x\oplus y = x^\sim \to y^{\sim-}= y^- \squig
x^{\sim-}$ that corresponds to an ``MV-addition". And for any pseudo
MV-algebra $A$ we know, \cite{Dvu2}, that an MV-state is a
state-morphism iff $m(a\oplus b) = m(a)\oplus_{\mbox{\tiny \L}}
m(b)$ for all $a,b \in A$. Inspired by this we can characterize
state-morphisms as follows.

\begin{lemma} \label{psBCK-2.330} Let $m$ be a Bosbach state  on a bounded
pseudo-BCK algebra $A.$ The
following statements are equivalent:\\
$(a)$  $m$ is a state-morphism.\\
$(b)$  $m(a^\sim \to b^{-\sim})= \min\{m(a)+m(b),1\}$
for all $a,b\in A$.\\
$(c)$    $m(b^- \squig a^{\sim-})= \min\{m(a)+m(b),1\}$ for all
$a,b\in A$.
\end{lemma}

\begin{proof}  Due to Theorem \ref{psBCK-2.310}, conditions $(b)$
and $(c)$ are equivalent, and, moreover, we can assume that $m$ is
the same as  $\hat m.$

$(a)\Rightarrow (b).$  Assume that $m$ is a state-morphism on $A$, so it is a Bosbach state.\\
By Proposition \ref{psBCK-2.275-20}(3), we have
\begin{eqnarray*}
m(a^{\sim}\rightarrow b^{- \sim}) &=&m(a^{\sim}\rightarrow
b)=m(a^{\sim})\rightarrow_{\mbox{\tiny \L}} m(b)\\&=&
m(a)^-\rightarrow_{\mbox{\tiny \L}} m(b)=\min\{m(a)+m(b),1\}.
\end{eqnarray*}

$(b)\Rightarrow (a).$  Exchanging $m$ with $\hat m,$ we have
\begin{eqnarray*}
m(a\rightarrow b)&=&m(a^{- \sim}\rightarrow b^{-
\sim})=\min\{m(a^-)+m(b^{\sim-}),1\}\\
&=& \min\{1-m(a)+m(b),1\}= m(a)\rightarrow_{\mbox{\tiny \L}} m(b).
\end{eqnarray*}
Similarly, $m(a\rightsquigarrow b)=m(a)\rightarrow_{\mbox{\tiny \L}}
m(b)$. Hence, $m$ is a state-morphism.
\end{proof}

\begin{prop} \label{psBCK-2.340}
Let $s$ be a Bosbach state on  a bounded pseudo-BCK algebra
$A.$ The following are equivalent:\\
$(a)$ $s$ is a state-morphism.\\
$(b)$ $\Ker(s)$ is a normal and maximal filter of $A$.
\end{prop}

\begin{proof} $(a)\Rightarrow (b).$ It is similar as in \cite[Prop.3.4.10]{Kuhr}.\\
$(b)\Rightarrow (a).$ Let $K=\Ker(s).$  According to Theorem
\ref{psBCK-2.310}, $A/K= (A/K)/\Ker(\hat s)$ is a BCK algebra that
is term equivalent to an MV-algebra. Assume $F$ is a filter of $A/K$
and let $K(F)=\{a\in A\mid  a/K \in F\}.$ Then $K(F)$ is a filter of
$A$ containing $K.$  The maximality of $K$ implies $K=K(F)$ and $F
=\{1/K\}.$  Due to Theorem \ref{psBCK-2.310}, $A/K$ can be assumed
to be an archimedean MV-algebra having only one maximal filter,
$\{1/F\}$. Therefore, $A/K$ is an MV-subalgebra of the MV-algebra of
the real interval $[0,1].$  This yields that the mapping $a\mapsto
a/K$ $(a\in A)$ is the Bosbach state $s$ that is a state-morphism.
\end{proof}

\begin{lemma} \label{psBCK-2.350} Let $m$ be a state-morphism on a bounded
pseudo-BCK algebra $A$ and $K=\Ker(m)$. Then \\
$\hspace*{4cm}$ $a/K\leq b/K$ if and only if $m(a)\leq m(b),$  \\
$\hspace*{4cm}$ $a/K=b/K$ if and only if $m(a)=m(b)$.
\end{lemma}

\begin{proof} By Proposition \ref{psBCK-2.280-20} it follows that $m$ is a Bosbach state on $A$.
Applying Lemma \ref{psBCK-2.300}, $a/K\leq b/K$ iff $m(b)=m(a\vee_1 b)$. \\
But $m(a\vee_1 b)=\max\{m(a),m(b)\}$ and hence $m(a)\leq m(b)$.\\
For the later assertion we apply the first one.
\end{proof}

\begin{prop} \label{psBCK-2.360}
Let $m$ be a state-morphism on a bounded pseudo-BCK algebra $A$.
Then $(m(A),\oplus, ^-, 0)$ is a subalgebra of the standard
MV-algebra $([0,1],\oplus,^-,0)$ and the mapping $a/\Ker(m) \mapsto
m(a)$ is an isomorphism of $A/\Ker(m)$ onto $m(A)$.
\end{prop}

\begin{proof}
Similarly as \cite[Prop. 3.4.12]{Kuhr}.
\end{proof}

\begin{prop} \label{psBCK-2.370} Let $A$
be a bounded pseudo-BCK algebra and $m_1,m_2$ two state-morphisms
such that $\Ker(m_1)=\Ker(m_2)$. Then $m_1=m_2$.
\end{prop}

\begin{proof}  By Proposition \ref{psBCK-2.280-20}, $m_1$ and $m_2$
are two Bosbach states.  The conditions yield $A/\Ker(m_1) =
A/\Ker(m_2),$ and as in the proof of Proposition \ref{psBCK-2.340},
we have that $A/\Ker(m_1)$ is in fact an MV-subalgebra of the
MV-algebra of the real interval $[0,1].$ But $\Ker(\hat
m_1)=\{1/K\}=\Ker(\hat m_2).$   Hence, by \cite[Prop. 4.5]{Dvu2},
$\hat m_1 = \hat m_2,$ consequently, $m_1=m_2.$
\end{proof}

Let $A$ be a bounded pseudo-BCK algebra. We say that a  Bosbach
state $s$ is \emph{extremal} if for any $0<\lambda<1$ and for any
two  Bosbach states $s_1, s_2$ on $A$, $s=\lambda
s_1+(1-\lambda)s_2$ implies $s_1=s_2.$

Summarizing previous characterizations of state-morphisms, we have
the following result.

\begin{theo} \label{psBCK-2.410}
Let $s$ be a Bosbach on a bounded pseudo-BCK algebra $A.$
Then the following are equivalent:\\
$(a)$ $s$ is an extremal  Bosbach state.\\
$(b)$ $s(x \vee_1 y)=\max\{s(x),s(y)\}$ for all $x,y \in A.$ \\
$(c)$ $s(x \vee_2 y)=\max\{s(x),s(y)\}$ for all $x,y \in A.$ \\
$(d)$ $s$ is a state-morphism.\\
$(e)$ $\Ker(s)$ is a maximal filter.
\end{theo}
\begin{proof}  The equivalence of $(b)$--$(e)$ was proved in
Propositions \ref{psBCK-2.280-30} and \ref{psBCK-2.340}.

$(d)\Rightarrow (a).$  Let $s= \lambda s_1 +(1-\lambda)s_2,$ where
$s_1,s_2$ are Bosbach states and $0< \lambda <1.$  Then $\Ker(s) =
\Ker(s_1)\cap \Ker(s_2)$ and the maximality of $\Ker(s)$ gives that
$\Ker(s_1)$ and $\Ker(s_2)$ are maximal and normal filters. $(e)$
yields that $s_1$ and $s_2$ are state-morphisms and Proposition
\ref{psBCK-2.370} entails $s_1 = s_2 = s.$

$(a)\Rightarrow (d).$ Let $s$ be an extremal state on $A$. Define
$\hat s$ by  Proposition \ref{psBCK-2.300} on $A/\Ker(s)$. We assert
that $\hat s$ is an extremal MV-state on the MV-algebra $A/\Ker(s)$.
Indeed, let $\hat m = \lambda \mu_1 + (1-\lambda) \mu_2,$ where
$0<\lambda<1$ and $\mu_1$ and $\mu_2$ are states on $A/\Ker(s).$
There exist two Bosbach states $s_1$ and $s_2$ on $A$ such that
$s_i(a) := \mu_i(a/\Ker(s))$, $a\in A$ for $i=1,2.$ Then $s =
\lambda s_1 + (1-\lambda) s_2$ which gives $s_1 = s_2 = s,$ so that
$\mu_1 = \mu_2 = \hat s.$

Since $A/\Ker(s)$ is in fact an MV-algebra, we conclude from \cite[
Thm 6.1.30]{DvPu} that $\hat s$ is a state-morphism on $A/\Ker(s)$.
Consequently, so is $s$ on $A.$
\end{proof}

\begin{rem} \label{psBCK-2.380}
In the case of pseudo-BL algebras and bounded non-commutative
R$\ell$-monoids it was proved that the existence of a state-morphism
is equivalent with the existence of a maximal filter which is normal
(see \cite{Geo1} and respectively \cite{DvRa3}). This result is
based on the fact that, if $A$ is one the above mentioned structures
and $H$ is a maximal and normal filter of $A$, then $A/H$ is an
MV-algebra.

In the case of pseudo-BCK algebras this result is not true as we can
see in the next example.
\end{rem}

\begin{ex} \label{psBCK-2.390} Consider $A=\{0,a,b,c,1\}$ with $0<a<b<c<1$ and
the operations $\rightarrow,\rightsquigarrow$ given by the following tables:
\[
\hspace{10mm}
\begin{array}{c|ccccc}
\rightarrow & 0 & a & b & c & 1 \\ \hline
0 & 1 & 1 & 1 & 1 & 1 \\
a & b & 1 & 1 & 1 & 1 \\
b & b & c & 1 & 1 & 1 \\
c & 0 & a & b & 1 & 1 \\
1 & 0 & a & b & c & 1
\end{array}
\hspace{10mm}
\begin{array}{c|ccccc}
\rightsquigarrow & 0 & a & b & c & 1 \\ \hline
0 & 1 & 1 & 1 & 1 & 1 \\
a & b & 1 & 1 & 1 & 1 \\
b & b & b & 1 & 1 & 1 \\
c & 0 & b & b & 1 & 1 \\
1 & 0 & a & b & c & 1
\end{array}
.
\]
Then $(A,\leq,\rightarrow,\rightsquigarrow)$ is a good pseudo-BCK
lattice. $D=\{1\}$ is a maximal normal filter of $A$, but there is
no Bosbach state on $A$.  Indeed, assume that $A$ admits a Bosbach
state $s$ such that $s(0)=0$, $s(a)=\alpha$, $s(b)=\beta$,
$s(c)=\gamma$, $s(1)=1$. From $s(x)+s(x\rightarrow
y)=s(y)+s(y\rightarrow x)$, taking $x=a, y=0$, $x=b,y=0$ and
respectively $x=c,y=0$ we get $\alpha=1/2$, $\beta =
1/2$, $\gamma = 1$. \\
On the other hand, we have:\\
$\hspace*{3cm}$ $s(a)+s(a\rightsquigarrow b)=s(a)+s(1)=1/2+1=3/2,$\\
$\hspace*{3cm}$ $s(b)+s(b\rightsquigarrow a)=s(b)+s(b)=1/2+1/2=1,$\\
so condition $(B_2)$ does not hold.\\
Thus, there is no Bosbach state, in particular, no state-morphism on
$A$.
\end{ex}

\begin{rem} \label{psBCK-2.395}
The reason why the above result  does not hold is that $A/\{1\}\cong
A$ is not an MV-algebra, see Remark \ref{psBCK-2.380}.
\end{rem}

Inspired by the latter remark, we have the following
characterization of the existence of Bosbach states on a bounded
pseudo-BCK algebra.

\begin{theo}\label{thAD}  Let $A$ be a bounded pseudo-BCK algebra.
The following statements are equivalent:\\
$(a)$ $A$ admits a Bosbach state.\\
$(b)$ There exists a normal filter $F\ne A$ of $A$ such that $A/F$
is termwise equivalent to an MV-algebra.\\
$(c)$  There exists a normal and maximal filter $F$ such that $A/F$
is termwise equivalent to an MV-algebra.
\end{theo}

\begin{proof} $(a)\Rightarrow (b).$  Let $m$ be a Bosbach state,
then the normal filter $F=\Ker(m)$ according to Theorem
\ref{psBCK-2.310} satisfies $(b).$

$(b)\Rightarrow (a).$ If $A/F$ is an MV-algebra, then it possesses
at least one MV-state, say $\mu.$ The function $m(a):=\mu(a/F)$
$(a\in A)$ is a Bosbach state on $A$.

$(a)\Rightarrow (c).$ If $A$ possesses at least one state, by the
Krein-Mil'man theorem (see (3.1) below), $\partial_e \mathcal{BS}(A)
= \mathcal{SM}(A) \ne \emptyset.$ Then there is a state-morphism $m$
on $A$ and due to Theorem \ref{psBCK-2.410}(d), the filter
$F=\Ker(m)$ is maximal and normal, and by Theorem \ref{psBCK-2.310},
$F$ satisfies $(c).$

$(c)\Rightarrow (a).$ It is the same as that of $(b)\Rightarrow
(a).$
\end{proof}

\begin{rem} \label{psBCK-2.400}
The previous example of a linearly ordered stateless pseudo-BCK
algebra shows another difference between pseudo-BCK algebras and
pseudo-BL algebras because: \emph{Every linearly ordered pseudo-BL
algebra admits a  Bosbach states} (see
\cite{Dvu3, DvHy2}). \\
\end{rem}

We say that a net of  Bosbach states $\{s_\alpha\}$ {\it converges
weakly} to a  Bosbach state $s$ if $s(a) = \lim_\alpha s_\alpha(a)$
for every $a \in A.$ According to the definition of Bosbach states,
the set of Bosbach states  is a compact Hausdorff topological space
(theoretically empty) in the weak topology.

Extremal Bosbach states are very important because they generate all
Bosbach states:  Due to the Krein--Mil'man theorem, \cite[Thm
5.17]{Goo}, every Bosbach state is a weak limit of  a net of convex
combinations of extremal Bosbach states.

Let $\mathcal {BS}(A),$ $\partial_e\mathcal {BS}(A)$ and $\mathcal
{SM}(A)$ denote the set of all Bosbach states, all extremal Bosbach
states, and all state-morphisms on $(A,\to, \rightsquigarrow, 0,1),$
respectively.  Theorem \ref{psBCK-2.410} says
$$ \partial_e\mathcal {BS}(A) = \mathcal {SM}(A)\eqno(3.1)
$$
and they are compact subsets of $\mathcal {BS}(A)$ in the weak
topology.

\begin{Def} \label{psBCK-2.430} $\rm($\cite{Ciu12}$\rm)$ Let $A$ be a good
bounded pseudo-BCK algebra. The elements $x,y \in A$ are called
\emph{orthogonal}, denoted by $x \perp y$ iff $x^{- \sim} \leq
y^{\sim}$. If the elements $x,y \in A$ are orthogonal, we define a
partial operation $+$ on $A$ by $x+y:=x \oplus y$.
\end{Def}

\begin{Def}\label{psBCK-2.450} $\rm($\cite{Ciu12}$\rm)$
Let $A$ be a good bounded pseudo-BCK algebra. A {\it Rie\v can
state} on $A$ is a function $s:A\longrightarrow [0,1]$
such that the following conditions hold for all $x,y\in A:$\\
$(R_1)$ If $x \perp y$, then $s(x+y)=s(x)+s(y);$\\
$(R_2)$ $s(1)=1.$
\end{Def}

\begin{prop} \label{psBCK-2.100} $\rm($\cite{Ciu12}$\rm)$
If $s$ is a Rie\v{c}an state on a good bounded pseudo-BCK algebra $A$, then the
following properties hold for all $x,y\in A:$\\
$(1)$ $s(x^{-})=s(x^{\sim})=1-s(x).$\\
$(2)$ $s(0)=0.$ \\
$(3)$ $s(x^{- \sim})=s(x^{\sim -})=s(x^{- -})=s(x^{\sim \sim})=s(x).$ \\
$(4)$ if $x\leq y$, then $s(x)\leq s(y)$ and
         $s(y \rightarrow x^{- \sim})=s(y \rightsquigarrow x^{\sim -})=1+s(x)-s(y).$\\
$(5)$ $s((x \vee_1 y) \rightarrow x^{- \sim})=s((x \vee_1 y)
\rightsquigarrow x^{- \sim})=1-s(x\vee_1 y)+s(x)$
and \\
$\hspace*{0.5cm}$
$s((x \vee_2 y) \rightarrow x^{- \sim})=s((x \vee_2 y) \rightsquigarrow x^{- \sim})=1-s(x\vee_2 y)+s(x).$\\
$(6)$ $s((x \vee_1 y) \rightarrow y^{- \sim})=s((x \vee_1 y) \rightsquigarrow y^{- \sim})=1-s(x\vee_1 y)+s(y)$ and \\
$\hspace*{0.5cm}$ $s((x \vee_2 y) \rightarrow y^{- \sim})=s((x
\vee_2 y) \rightsquigarrow y^{- \sim})=1-s(x\vee_2 y)+s(y)$.
\end{prop}

\begin{rem} According to \cite[Thm 3.17]{Ciu12}  every Bosbach state on a good bounded pseudo-BCK algebra
is a Rie\v{c}an state.  The converse is not true in general, as it
was proved in \cite[Ex. 3.18]{Ciu12}, see also the next Example
\ref{psBCK-2.120}.
\end{rem}

In the following example, we show that there exists a Rie\v can
state on a bounded pseudo-BCK algebra $A,$ that is not a Bosbach
state. Moreover, this  $A$ has no Bosbach state.

\begin{ex} \label{psBCK-2.120}
Consider $A=\{0,a,b,c,1\}$ with $0<a<b<c<1$  from Example
\ref{psBCK-2.390}. Then $(A,\leq,\rightarrow,\rightsquigarrow)$ is a
good pseudo-BCK algebra.
The function $s : A \longrightarrow [0,1]$ defined by \\
$\hspace*{3cm}$ $s(0)=0, s(a)=1/2, s(b)=1/2, s(c)=1, s(1)=1$ \\
is a unique Rie\v can state. Indeed, the orthogonal elements of $A$
are the pairs $(x,y)$ in the following table:
$$\begin{tabular}{|l|l|l|l|l|}
\hline $x$ & $y$ &  $x^{- \sim}$ & $y^{\sim}$ & $x+y$\\
\hline $0$ & $0$ & $0$ & $1$ & $0$\\
\hline $0$ & $a$ & $0$ & $b$ & $b$\\
\hline $0$ & $b$ & $0$ & $b$ & $b$\\
\hline $0$ & $c$ & $0$ & $0$ & $1$\\
\hline $0$ & $1$ & $0$ & $0$ & $1$\\
\hline $a$ & $0$ & $b$ & $1$ & $b$\\
\hline $a$ & $a$ & $b$ & $b$ & $1$\\
\hline $a$ & $b$ & $b$ & $b$ & $1$\\
\hline $b$ & $0$ & $b$ & $1$ & $b$\\
\hline $b$ & $a$ & $b$ & $b$ & $1$\\
\hline $b$ & $b$ & $b$ & $b$ & $1$\\
\hline $c$ & $0$ & $1$ & $1$ & $1$\\
\hline $1$ & $0$ & $1$ & $1$ & $1$\\
\hline\end{tabular}$$
We prove now that $s$ is not a Bosbach state. Look at $a<b$. We
have: $s(b\rightarrow a)=s(c)=1,$ but $s(b\rightsquigarrow
a)=s(b)=1/2.$ So, $s$ is not a Bosbach state. We recall that the
kernel of $s$ is the set $\{c,1\}$ that is not a filter.

On the other hand, as it was shown in Examle \ref{psBCK-2.390}, $A$
has no Bosbach state.
\end{ex}

\section{Measures on Pseudo-BCK Algebras}

In this section we generalize  measures on pseudo-BCK algebras
introduced by A. Dvure\-\v{c}en\-skij in \cite{Dvu1} and \cite{DvPu}
to pseudo-BCK algebras that are not necessarily bounded. In
particular, we show that if $A$ is a downwards-directed pseudo-BCK
algebra and $m$ a measure on it, then the quotient over the kernel
of $m$ can be embedded into the negative cone of an abelian,
archimedean $\ell$-group as its subalgebra.  This result will enable
us to characterize nonzero measure-morphisms as measures whose
kernel is a maximal filter.

Consider the bounded $\vee_1$-commutative BCK(P) algebra (i.e. an
MV-algebra) ${\mathcal A}_{\mbox{\tiny \L}}=([0,1],\leq,
\rightarrow_{\mbox{\tiny \L}},0,1)$, where $\rightarrow_{\mbox{\tiny
\L}}$ is the \L ukasiewicz implication: $x \rightarrow_{\mbox{\tiny
\L}} y = \min\{1-x+y,1\}$.

\begin{Def}\label{psBCK-3.10} Let $(A,\leq,\rightarrow,
\rightsquigarrow,1)$ be a pseudo-BCK algebra.
A mapping $m:A\longrightarrow [0, \infty)$ such that for all $x,y \in A$, \\
$(1)$ $m(x\rightarrow y)=m(x\rightsquigarrow y)=m(y)-m(x)$ whenever
$y\leq x$ is said to be a \emph{measure};\\
$(2)$ if $0\in A$  and $m$ is a measure with $m(0)=1$, then $m$ is
said
to be a \emph{state-measure};\\
$(3)$ if $m(x\rightarrow y)=m(x\rightsquigarrow
y)=\max\{0,m(y)-m(x)\}$ is said to be a \emph{measure-morphism};\\
$(4)$  if $0\in A,$ $m(0)=1$ and  $m$ is a measure-morphism, $m$ is
said to be a {\it state-measure-morphism}.
\end{Def}

Of course, the function vanishing on $A$ is always a (trivial)
measure.

We note that our definition of a measure (a state-measure) is a
definition of a map that maps pseudo-BCK algebra that is in the
``negative cone" to the positive cone of the reals $\mathbb R.$  For
a relationship with the previous type of Bosbach state see the
second part of the present section and Remark \ref{psBCK-3.90}.

For example, let $(G,+,0,\wedge,\vee)$ be an $\ell$-group, with the
negative cone $G^-,$ see Example \ref{psBCK-1.10-20}. Assume that
$m$ is a positive-valued function on $G^-$ that preserves addition
in $G^-.$  Then $m$ is a measure on the pseudo-BCK algebra $G^-,$
and conversely if $m$ is a measure  on $G^-,$  then $m$ is additive
on $G^-$ and positive-valued.

We recall that not every negative cone even  of abelian $\ell$-group
admits a nontrivial measure. To see that look  \cite[Ex. 9.6]{Goo}.

\begin{prop}\label{AD5} Let $m$ be a measure on a pseudo-BCK algebra
$A$.   For all $x,y \in A$, we have\\
$(1)$ $m(1)=0.$ \\
$(2)$ $m(x)\geq m(y)$ whenever $x\leq y.$  \\
$(3)$ $m(x\vee_1 y) = m(y\vee_1 x)$ and $m(x\vee_2 y) = m(y\vee_2
x).$\\
$(4)$ $m(x\vee_1 y) = m(x\vee_2 y).$\\
$(5)$ $m(x\to y) = m(x\squig y).$
\end{prop}

\begin{proof}
(1) Since $1\leq 1$ we get $m(1)=m(1\rightarrow 1)=m(1)-m(1)=0.$

(2) Since $x\leq y$ it follows that $m(y\rightarrow x)=m(x)-m(y)$,
so $m(x)-m(y)\geq 0.$

(3) First, let $ x\le y$. Then by Proposition
\ref{psBCK-1.20-10}(2), we have $m(x\vee_1 y) = m(y).$ Using the
property of measures, we have $m((y\vee_1 x)\to x) = m(x)-m(y\vee_1
x) = m(x) - m((y\to x)\squig x) = m(x) -m(x) +m(y\to x)= m(x)-m(y)$
giving $m(y\vee_1 x) = m(y).$

Let now $x,y\in A$ be arbitrary. Using the first part of the present
proof and (2), we have $m(x\vee_1 y) = m(x\vee_1 (x \vee_1 y)) =
m((x\vee_1 y) \vee_1 x) \le m(y \vee_1 x) = m(y \vee_1 (y\vee_1 x))
= m((y\vee_1 x)\vee_1 y) \le m(x\vee_1 y).$

In a similar way we prove the second equation in (3).

(4) First again, let $x \le y$.  Then $m(x \vee_1 y) = m(y).$ And
$m(x\vee_2 y) = m(y\vee_2 x) = m(y).$

Let now $x,y\in A$ be arbitrary. Using (3), we have: $m(x\vee_1 y) =
m(x\vee_1 (x\vee_1 y)) = m((x\vee_1 y)\vee_1 x) = m((x\vee_1
y)\vee_2 x) \le m(y\vee_2 x) = m(x\vee_2 y) = m(x\vee_2 (x\vee_2 y))
= m((x\vee_2 y) \vee_2 x)= m((x\vee_2 y)\vee_1 y) \le m(x\vee_1 y).$

$(5)$ According to (2.1) and (4), $m(x \to y) = m((x\vee_1 y)\to y)
= m(y) - m(x\vee_1 y)= m(y) - m(x\vee_2 y) = m((x\vee_2y)\squig y) =
m(x\squig y).$
\end{proof}

\begin{prop}\label{psBCK-3.20} Let $A$ be a
pseudo-BCK algebra. Then  \\
$(1)$ $y\leq x$ implies $m((x\rightarrow y)\rightsquigarrow
y)=m((x\rightsquigarrow y)\rightarrow y)=m(x)$ whenever $m$ is a
measure on $A.$ \\
$(2)$ If $m$ is a measure on $A$, then $\Ker_0(m)=\{x\in A \mid
m(x)=0\}$ is a normal filter  of $A.$ \\
$(3)$ Any measure-morphism on $A$ is a measure on $A.$
\end{prop}
\begin{proof}
$(1)$ From $y\leq x\rightarrow y$  we get
      $m((x\rightarrow y)\rightsquigarrow y)=m(y)-m(x\rightarrow y)=m(y)-(m(y)-m(x))=m(x).$\\
      Similarly, $m((y\rightsquigarrow x)\rightarrow x)=m(x).$ \\
$(2)$ According to Proposition \ref{AD5}(1), $1\in \Ker_0(m)$. Assume that $x,x\rightarrow y \in \Ker_0(m)$. \\
      Since $x\leq x \vee_1 y$, due to Proposition \ref{AD5}(2), we have $0=m(x)\geq m(x\vee_1 y)$, so $m(x\vee_1 y)=0$.\\
      In addition,\\
      $0=m(x\rightarrow y)=m(x\vee_1 y\rightarrow y)=m(y)-m(x\vee_1 y)=m(y)$, so $y\in
      \Ker_0(m)$.\\
      (we applied the fact that $y\leq x\vee_1 y$ and Proposition \ref{psBCK-1.20-20}). \\
      Thus, $\Ker_0(m)$ is a filter of $A.$ The normality of $\Ker_0(m)$ follows from
      Proposition \ref{AD5}(5).\\
$(3)$ We have $m(1)=m(1\to 1) = \max\{0, m(1)-m(1)\}=0$ so that if
$y\le x$ then $0=m(1)=m(y\to x)= \max\{0, m(x)-m(y)\}$ and $m(x)\le
m(y)$ so that $m(x\to y) = \max\{0, m(y)-m(x)\}= m(y)-m(x).$
Similarly, $m(x\rightsquigarrow y) = m(y)-m(x).$
\end{proof}

\begin{ex} \label{psBCK-3.25} Consider the bounded pseudo-BCK lattice $A$ from Example \ref{psBCK-1.190}.\
The function $m : A \longrightarrow [0,\infty)$ defined by: $m(0)=1,
m(a)=m(b)=m(c)=m(d)=m(1)=0$ is a unique measure on $A$. Moreover,
$m$ is even a state-measure on $A$.
\end{ex}

\begin{prop} \label{psBCK-3.60} Let $A$
be a bounded pseudo-BCK algebra. If $M$ is a Bosbach state, then
$m=1-M$ is a state-measure.
\end{prop}

\begin{proof}
Let $y\leq x$, that is $y\rightarrow x=y\rightsquigarrow x=1$.\\
Replacing in the axioms $(B1)$ and $(B2)$ of a Bosbach state, we obtain:\\
$M(x)+M(x\rightarrow y)=M(y)+1$ and
$M(x)+M(x\rightsquigarrow y)=M(y)+1,$\\
which implies $M(x\rightarrow y)=M(x\rightsquigarrow y)=1-M(x)+M(y).$\\
We get $m(x\rightarrow y)=m(x\rightsquigarrow y)=1-M(x\rightarrow y)=M(x)-M(y)=1-M(y)-(1-M(x))=m(y)-m(x)$.\\
Also, $m(0)=1-M(0)=1$ by $(B3),$ and thus, $m$ is a state-measure.
\end{proof}

\begin{prop} \label{psBCK-3.60-10}
Let $A$ be a bounded pseudo-BCK algebra.
 If $m$ is a state-measure on
$A$, then $M=1-m$ is a Bosbach state on $A.$
\end{prop}

\begin{proof}
We have: $y\leq x\vee_1 y$ and using the definition of the measure, we get\\
$m(x\vee_1 y \rightarrow y)=m(x\vee_1 y \rightsquigarrow
y)=m(y)-m(x\vee_1 y)$. Using Proposition \ref{psBCK-1.20-20}, we
have: $x\vee_1 y\rightarrow y=x\rightarrow y$, so we get
$m(x\rightarrow y)=m(y)-m(x\vee_1 y).$\\ Similarly, $m(y\rightarrow
x)=m(x)-m(y\vee_1 x)$ and also,
$m(x\rightsquigarrow y)=m(y)-m(x\vee_2 y),$\\
$m(y\rightsquigarrow x)=m(x)-m(y\vee_2 x).$ According to Proposition
\ref{AD5}(3), $m(x\vee_1 y)=m(y\vee_1 x),$ so that
$m(x)+m(x\rightarrow y)=m(y)+m(y\rightarrow x)$.\\
Similarly, $m(x)+m(x\rightsquigarrow y)=m(y)+m(y\rightsquigarrow x)$.\\
Therefore, $M(x)+M(x\rightarrow y)=M(y)+M(y\rightarrow y)$ and
$M(x)+M(x\rightsquigarrow y)=M(y)+M(y\rightsquigarrow y)$.
Also, $M(0)=0$ by the hypothesis, and $M(1)=1$ by Proposition \ref{AD5}(1).\\
Thus, $M$ is a Bosbach state.
\end{proof}

If $A$ is a bounded pseudo-BCK algebra, in a similar way as for
Bosbach states, we can define extremal state-measures, as well as
the weak-topology. Let us denote the set of state-measures,
$\mathcal{SM}_1(A),$ the set of state-measure-morphisms,
$\mathcal{SMM}_1(A),$ and the set of extremal state-measures,
$\partial_e \mathcal{SM}_1(A),$ respectively.

\begin{theo}\label{AD6}  Let $A$ be a bounded
pseudo-BCK algebra. Define a map $\Psi:\ \mathcal{SM}_1(A) \to
\mathcal{BS}(A)$ via $\Psi(m) = 1-m,$ $m\in \mathcal{SM}_1(A).$ Then
$\Psi$ is a affine-homeomorhism such that $m$ is a
state-measure-morphism if and only if $\Psi(m)$ is a state-morphism.
In particular, $m$ is an extremal state-measure if and only if $m$
is a state-measure-morphism.
\end{theo}

\begin{proof}  Propositions \ref{psBCK-3.60}-\ref{psBCK-3.60-10}
show that $\Psi$ is a bijection preserving convex combinations and
weak topologies.

If, say, $s$ is a state-morphism on $A,$ i.e. $s(x\to y)
=\min\{1-m(x)+m(y)\},$ then it is straightforward to show that
$m=1-s$ is a state-measure-morphism on $A,$ i.e. $m(x\to
y)=\max\{m(y)-m(x),0\}$ (as well as for the second arrow $\squig$).

In view of Theorem \ref{psBCK-2.410}, we see that a state-measure is
extremal iff it is state-measure-morphism.
\end{proof}

As a corollary of Theorem \ref{AD6} and (3.1), we have that if $A$
is a bounded pseudo-BCK algebra, then
$$
\partial_e\mathcal{SM}_1(A) = \mathcal{SMM}_1(A). \eqno(4.1)
$$

\begin{theo} \label{psBCK-3.50}  Let $m$ be a measure on a
pseudo-BCK algebra $A.$  Then $A/\Ker_0(m)$ is a pseudo-BCK algebra
and the mapping $\hat{m}:A/\Ker_0(m) \longrightarrow [0,+\infty)$
defined by $\hat{m}(\bar{x}):=m(x),$ $\bar{x}:=x/\Ker_0(m)\in
A/\Ker_0(m),$ is a measure on $A/\Ker_0(m),$ and $A/\Ker_0(m)$ is
$\vee$-commutative.
\end{theo}

\begin{proof}
By (2.1) and (4) of Proposition \ref{AD5}, we have $m(x\to y) =
m(x\vee_1 y \to y)=m(y)-m(x\vee_1 y) = m(y)-m(x\vee_2 y) = m(x\vee_2
y \rightsquigarrow y) = m(x\rightsquigarrow y).$ According to (2) of
Proposition \ref{psBCK-3.20}, $\Ker_0(m)$ is a normal filter. Since
if $x\to y, y\to x \in \Ker_0(m),$ we have  $m(x\to y) = m(x\vee_1
y\to y) = 0 = m(y) -m(x\vee_1y)$. Similarly, $m(x)=m(y\vee_1 x).$
But $m(y)=m(x\vee_1 y) = m(y\vee_1 x) = m(x).$ Hence, $\hat m$ is a
well-defined function on $A/\Ker_0(m).$

We recall that $\bar x= \bar y$ iff $m(x)=m(y)= m(x\vee_1 y),$ and
$\bar x \le \bar y$ iff $m(x\to y)=0.$ Therefore, if
$x/\Ker_0(m)=y/\Ker_0(m),$ then $x/\Ker_0(m)= (x\vee_1 y)/\Ker_0(m)
=(x\vee_2 y)/\Ker_0(m)= (y\vee_1 x)/\Ker_0(m)=(y\vee_2
x)/\Ker_0(m).$

To show that $\hat m$ is a measure, assume $\bar y \le \bar x$. By
(2) of Proposition \ref{psBCK-1.20-10}, $\bar y \vee_1 \bar x= \bar
x.$ Then $\hat m (\bar x \to \bar y)= m(x\to y) = m(x\vee_1 y\to y)
= m(y) - m(x\vee_1 y).$ But $m(x\vee_1 y) = m(y \vee_1 x) = \hat
m(\bar y\vee_1 \bar x) = \hat m(\bar x) = m(x).$ Therefore, $\hat
m(\bar x \to \bar y) = \hat m(\bar y)-\hat m(\bar x).$ Similarly,
$\hat m(\bar x \rightsquigarrow \bar y)=\hat m(\bar y)-\hat m(\bar
x).$

In the same way as in the proof of Proposition \ref{AD2} we can show
that $A/\Ker_0(M)$ is  $\vee$-commutative.
\end{proof}

In view of Theorem \ref{psBCK-2.310} and Theorem \ref{AD6} we know
that if $m$ is a state-measure on a bounded pseudo-BCK algebra $A,$
then $A/\Ker(m)$ is in fact an MV-algebra, so that according to the
famous representation theorem of Mundici, \cite{CDM}, $A/\Ker(m)$ is
an interval in an $\ell$-group with strong unit. In the following
result we generalize this $\ell$-group representation of the
quotient for measures on  unbounded pseudo-BCK algebra that are
downwards-directed.

\begin{theo}\label{ADThm1} Let $m$ be a measure on an unbounded pseudo-BCK algebra $A$ that is
a downwards-directed set. Then the arrows $\to/\Ker_0(m)$ and $
\rightsquigarrow/\Ker_0(m)$ on $A/\Ker_0(m)$ coincide. Moreover,
there is a unique (up to isomorphism)) archimedean $\ell$-group $G$
such that $A/\Ker_0(m)$ is a subalgebra of the pseudo-BCK algebra
$G^-$ and $A/\Ker_0(m)$ generates the $\ell$-group $G.$
\end{theo}

\begin{proof}  We note that if $a$ is an arbitrary element of $A$,
then $([a,1],\le, \to, \squig,a,1)$ is a pseudo-BCK algebra.

We denote by $K_0:=\Ker_0(m).$  Given $x,y\in A,$ choose an element
$a \in A$ such that $a\le x,y.$

If $m(a)=0$, then $a/K_0 = x\to y/K_0 = x\squig y/K_0 = 1/K_0.$

Assume $m(a) >0$ and define $m_a(z):= m(z)/m(a)$ for any $z \in
[a,1].$ Then $m_a$ is a state-measure on $[a,1]$ and in view of
Theorem \ref{AD6}, $s_a:=1-m_a$ is a Bosbach state on $[a,1].$
Theorem \ref{psBCK-2.310} entails that $[a,1]/\Ker(s_a)$ can be
converted into an archimedean MV-algebra.  In particular, $(x\to
y)/\Ker(s_a) = (x\squig y)/\Ker_0(s_a).$ This yields $s_a((x\to
y)\to (x\squig y))=1$ and $m((x\to y)\to (x\squig y))=0.$ In a
similar way, $m((x\squig y)\to (x\to y))=0.$ This proves
$\to/\Ker_0(m)=\squig /\Ker_0(m).$

In  addition, we can prove that, for all $x,y \in A,$
$$
((x\to y)\vee (y\to x))/\Ker_0(m) = 1/\Ker_0(m) = ((x\squig y)\vee (y\squig
x))/\Ker_0(m).
$$

It is clear that if $m=0$, then $\Ker_0(m)=A$ and
$A/\Ker_0(m)=\{1/\Ker_0(m)\}$ so that the trivial $\ell$-group
$G=\{0_G\},$ where $0_G$ is a neutral element of $G,$ satisfies our
conditions.

Therefore, let $m\ne 0.$ By  \cite[Lem 4.1.8]{Kuhr}, $A/\Ker_0(m)$
is a distributive lattice. As in Proposition \ref{psBCK-1.300} we
can show that $A/\Ker_0(m)$ satisfies the (RCP) condition, and
therefore, $A/\Ker_0(m)$ is a \L ukasiewicz BCK algebra, see
\cite{DvVe}. Therefore, \cite{DvVe, DvGr}, there is a unique (up to
isomorphism of $\ell$-groups) $\ell$-group $G$ such that
$A/\Ker_0(m)$ can be embedded into the pseudo-BCK algebra of the
negative cone $G^-,$ moreover, $A/\Ker_0(m)$ generates $G$. Since
the arrows in $A/\Ker_0(m)$ coincide, we see that $G$ is abelian,
and every interval $[a/K_0,1/K_0]$ is an archimedean MV-algebra, so
is $G.$
\end{proof}

We note that if $m$ is a measure-morphism on $A,$  then
$$
m(u\to^n x) = \max\{0,m(x)-nm(u)\}\eqno(4.2)
$$
for any $n\ge 0,$ and
$$
m(x_1\to(\cdots \to (x_n\to a)\cdots ))=\max\{0,
m(a)-m(x_1)-\cdots-m(x_n)\}.\eqno(4.3)
$$

\begin{prop}\label{pr:6.6}   Let $m$ be a measure-morphism on a
pseudo-BCK algebra $A$ such that $m\ne 0.$  Then $\Ker_0(m)$ is a
normal and maximal filter of $A.$
\end{prop}

\begin{proof} Since $m$ is a measure-morphism, by Proposition \ref{psBCK-3.20},
$\Ker_0(m)$ is a normal filter.

Choose $a \in A$ such that $m(a)\ne 0.$ Let $F$ be the filter
generated by $\Ker_0(m)$ and by the element $a.$ Let $z \in A$ be an
arbitrary element of $A.$ There is an integer $n\ge 1$ such that
$(n-1)m(a)\le m(z)< nm(a).$ Due to (4.3), we have that
$m(a\to^nz)=\max\{0,m(z)-nm(a)\}=0$ so that $z\in F$ and $A
\subseteq F$ proving that $\Ker_0(m)$ is a maximal filter.
\end{proof}

If $m\ne 0$ is a measure on a bounded pseudo BCK-algebra $A,$ then
passing to a state-measure $s_m(a):=m(a)/m(1),$ $a\in A,$ and using
Theorem \ref{AD6}, we see that $m$ is a measure-morphism iff
$\Ker_0(m)$ is a maximal filter. The same result is true for
unbounded pseudo-BCK algebra that is downwards-directed:

\begin{theo}\label{ADThm2}  Let $m\ne 0$ be a measure on an unbounded
pseudo-BCK algebra $A$ that is downwards directed. Then $m$ is a
measure-morphism if and only if $\Ker_0(m)$ is a maximal filter.
\end{theo}

\begin{proof}  By Proposition \ref{pr:6.6},  $\Ker_0(m)$ is a
maximal filter of $A.$

Suppose now $\Ker_0(m)$ is a maximal filter of $A.$ In view of
Theorem \ref{AD6}, $A/\Ker_0(m)$ can be embedded as a subalgebra
into the pseudo-BCK algebra $G^-$, where $G^-$ is the negative cone
of an abelian and archimedean $\ell$-group $G$ that is generated by
$A/\Ker_0(m).$ Let $\hat m(a/\Ker_0(m)):=m(a)$ $(a\in A).$ Then
$\Ker_0(\hat m) = \{1/\Ker_0(m)\}$ and $0_G:=1/\Ker_0(m)$ is the
neutral element of $G$.

Fix an element $a \in A$ with $m(a)>0.$ Since $\Ker_0(m)$ is maximal
in $A,$ $\Ker_0(\hat m)=\{1/\Ker_0(m)\}$ is maximal in $A/\Ker_0(m)$
and consequently, $\{1/\Ker_0(m)\}$ is a maximal filter of the
pseudo-BCK algebra $G^-$ because $A/\Ker_0(m)$ generates $G.$
Therefore, the $\ell$-ideal $L:=\{0_g\}=\{1/\Ker_0(m)\}$ is a
maximal $\ell$-ideal of $G.$ We recall that every maximal
$\ell$-ideal, $L,$ of an $\ell$-group is prime ($a,b \in G^+$ with
$a\wedge b = 0$ implies $a\in L$ or $b\in L$), whence $G/L$ is a
linearly ordered $\ell$-group (see e.g. \cite[Prop. 9.9]{Dar}).
Since $G = G/L,$ $G$ is archimedean and linearly ordered, due to the
H\"older theorem, \cite[Thm 24.16]{Dar}, $G$ is an $\ell$-subgroup
of the $\ell$-group of real numbers, $\mathbb R.$  Let $s$ be the
unique extension of $\hat m$ onto $G,$ then $s$ is additive on $G$
and $s(g)\ge 0$ for any $g \in G^-.$ Since $G$ is an $\ell$-subgroup
of $\mathbb R,$ $s$ is a unique additive function on $G$ that is
positive on the negative cone  (see example just after Definition
\ref{psBCK-3.10}) with the property $s(a/\Ker_0(m))=m(a)>0$ for our
fixed element $a \in A.$ Because $A/\Ker_0(m)$ can be embedded into
$\mathbb R^-,$ we see that $s$ is a measure-morphism on $G^-.$
Consequently, $m$ is a measure-morphism on $A.$
\end{proof}

\begin{prop}\label{ADProp1}  Let $m_1$ and $m_2$ be two
measure-morphisms  on a downwards-directed pseudo-BCK algebra $A$
such that there is an element $a\in A$ with $m_1(a)=m_2(a)>0.$ If
$\Ker_0(m_1)=\Ker_0(m_2),$ then $m_1 = m_2.$

In addition, let $a\in A$ be fixed. If $m$ is a measure-morphism on
$A$ such that $m(a)>0,$ then $m$ cannot be expressed as a convex
combination of two measure $m_1$ and $m_2$ such that
$m_1(a)=m_2(a)=m(a).$

\end{prop}

\begin{proof} (1) Due to Theorem \ref{ADThm2},
$A/\Ker_0(m_1)=A/\Ker_0(m_2)$ is a pseudo-BCK subalgebra of $\mathbb
R^-.$ The condition $m_1(a)=m_2(a)>0$ entails $\hat m_1 = \hat m_2$
so that $m_1 = m_2.$

(2) Let  $m=\lambda m_1 +(1-\lambda)m_2$ where $m_1$ and $m_2$ are
measures on $A$ such that $m_1(a)=m_2(a)=m(a)$ and $0<\lambda <1.$
Then $\Ker_0(m)\subseteq \Ker_0(m_1)\cap \Ker_0(m_2).$ The
maximality of $\Ker_0(m)$ entails that both $\Ker_0(m_1)$ and
$\Ker_2(m)$ are maximal ideals and by Theorem \ref{ADThm2}, we see
that $m_1$ and $m_2$ are measure-morphisms on $A.$  The condition
$m_1(a)=m_2(a)=m(a)$ yields by (1) that $m=m_1=m_2.$
\end{proof}

\begin{prop}  \label{psBCK-3.60-30}
Let $m$ be a state-measure on a good bounded pseudo-BCK algebra.\\
Then $M=1-m$ is a Rie\v can state.
\end{prop}

\begin{proof}
Let $x, y$ be a pair of orthogonal elements, that is $y^{- \sim}\leq
x^{-}$ and using the fact that $m$ is a measure, we obtain:
$m(x^{-}\rightarrow y^{- \sim})=m(x^{-}\rightsquigarrow y^{-
\sim})=m(y^{- \sim})-m(x^{-}).$\\ Now, because $A$ is good we get:
$m(x^{-}\rightarrow y^{\sim -})=m(y)-1+m(x),$ which implies
$M(x\oplus y)=M(x)+M(y).$  Thus, $M$ is a Rie\v can state.
\end{proof}

\begin{prop} \label{psBCK-3.70}
Let $A$ be a bounded pseudo-BCK(pDN) algebra and $s$ a Rie\v can
state on $A$. Then $S=1-s$ is a state-measure.
\end{prop}
\begin{proof} Let $s$ be a Rie\v can state on $A$. \\
Consider $y\leq x$. According to \cite{Ciu12} we have \\
$\hspace*{2cm}$
$s(x \rightarrow y^{- \sim})=s(x \rightsquigarrow y^{\sim -})=1-s(x)+s(y)$.\\
Taking into consideration the (pDN) condition we get \\
$\hspace*{2cm}$
$s(x \rightarrow y)=s(x \rightsquigarrow y)=1-s(x)+s(y)$.\\
It follows that $S(x \rightarrow y)=S(x \rightsquigarrow y)=S(y)-S(x)$.\\
Moreover, we have $S(0)=1$, so $S$ is a state-measure on $A$.
\end{proof}

\begin{rem} \label{psBCK-3.80}
We can also define a measure as a map $m:A\longrightarrow
(-\infty,0]$ such that \\
$\hspace*{2cm}$ $m(x\rightarrow y)=m(x\rightsquigarrow y)=m(x)-m(y)$
whenever $y\leq x$.\\
Properties (2) of Proposition \ref{AD5} and (1) in Proposition \ref{psBCK-3.20} become:\\
$(2^{\prime})$ $m(x)\leq m(y)$ whenever $x\leq y$ and $m$ is a measure on $A;$ \\
$(1^{\prime})$ $y\leq x$ implies $m((x\rightarrow y)\rightsquigarrow
y)=m((x\rightsquigarrow y)\rightarrow y)=-m(x)$ whenever $m$ is a
measure on $A;$ \\
If $m(0)=0$ then $m$ is a state on $A$.\\
Proposition \ref{psBCK-3.60} will be modified such that $m=1+M$.\\
Consider again the bounded pseudo-BCK lattice $A$ from Example
\ref{psBCK-1.190}.\ The function $m : A \longrightarrow (-\infty,0]$
defined by: $m(0)=-1, m(a)=m(b)=m(c)=m(d)=m(1)=0$ is the unique
measure on $A$.
\end{rem}

\begin{rem} \label{psBCK-3.90}
If a pseudo-BCK algebra is defined on the negative cone, like in
Examples \ref{psBCK-1.10-20} and \ref{psBCK-1.10-30}, we map through
the negative cone to the positive cone in ${\mathbb R}$. According
to the second kind of definition, we map the negative cone to
negative numbers.
\end{rem}

\section{Pseudo-BCK Algebras with Strong Unit}

In the present section, we will study state-measures on pseudo-BCK
algebras with strong unit. We apply the results of the previous
section to show how to characterize state-measure-morphisms as
extremal state-measures or as those with the maximal filter. In
particular, we show that for unital pseudo-BCK algebras that are
downwards directed, the quotient over the kernel can be embedded
into the negative cone of an abelian, archimedean $\ell$-group with
strong unit.

According to \cite{DvPu}, we are saying that an element $u$ of a
pseudo-BCK algebra $A$ is a {\it strong unit} if, for the filter
$(F9u)$ of $A$ that is generated by $u,$ we have $F(u)=A.$ For
example, if $(A,\to, \rightsquigarrow, 0,1)$ is a bounded pseudo-BCK
algebra, then $u=0$ is a strong element. If $G$ is an $\ell$-group
with strong unit $u\ge 0,$ then the negative cone $G^-$ is an
unbounded pseudo-BCK algebra with strong unit $-u.$

\begin{rem}\label{re:6.1} We note that a filter $F$ of a pseudo-BCK algebra with a strong
unit $u$ is a proper subset of $A$ if and only if $u \notin F.$
\end{rem}

By a {\it unital pseudo-BCK algebra} we mean a couple $(A,u)$ where
$A$ is a pseudo-BCK algebra with a fixed strong unit $u.$   We say
that a measure $m$ on $(A,u)$ is a {\it state-measure} if $m(u)=1.$
If, in addition, $m$ is a measure-morphism such that $m(u)=1$, we
call it also a {\it state-measure-morphism}.  We denote by $\mathcal
{SM}(A,u)$ and $\mathcal {SMM}(A,u)$ the set of all state-measures
and state-measure-morphisms on $(A,u),$ respectively. The set of
$\mathcal{SM}(A,u)$ is convex, i.e. if $m_1,m_2 \in
\mathcal{SM}(A,u)$ and $\lambda \in [0,1],$ then $m=\lambda
m_1+(1-\lambda)m_2 \in \mathcal{SM}(A,u);$ it could be empty.  A
state-measure $m$ is {\it extremal} if $m=\lambda
m_1+(1-\lambda)m_2$ for $\lambda \in (0,1)$ yields $m= m_1 = m_2.$
We denote by $\partial_e \mathcal{SM}(A,u)$ the set of all extremal
state-measures on $(A,u).$

\begin{ex}\label{ex:6.2}  Let $G$ be an $\ell$-group with  strong unit
$u\ge 0,$ i.e., given $g\in G$ there is an integer $n\ge 1$ such
that $g\le nu.$ Then a mapping $m$ on $G^-$ is a state-measure on
$(G^-,-u)$ if and only if (i) $m:\ G^- \to [0,\infty),$ (ii)
$m(g+h)= m(g)+m(h)$ for $g,h\in G^-,$ and (iii) $m(-u)=1.$ A
state-measure $m$ is extremal if and only if $m(g\wedge h)=
\max\{m(g),m(h)\},$ $g,h \in G^-,$ see \cite[Prop. 4.7]{Dvu2}. In
addition, $(-u)\to^n g = (g+nu)\wedge 0$ for any $n \ge 1.$
\end{ex}

\begin{ex}\label{ex:6.3}  Let $\Omega\ne \emptyset$ be a compact Hausdorff topological
space and let $\mbox{C}(\Omega)$ be the set of all continuous
functions on $\Omega.$ Then $\mbox{C}(\Omega)$ is an $\ell$-group
with respect to the pointwise ordering and usual addition of
functions and the element $u=1,$ the constant function equals $1$,
is a strong unit. According to  Riesz Representation Theorem, see
e.g. \cite[p. 87]{Goo},  a mapping $m:A \to [0,\infty)$ is a
state-measure on $(\mbox{C}(\Omega)^-,-1)$ if and only if there is a
Borel probability measure $\mu$ on $\mathcal B(\Omega)$ such that
$$
m(f)=-\int_\Omega f(x)\dx\mu(x),\quad f \in
\mbox{C}(\Omega)^-,\eqno(5.1)
$$
and vice-versa, given a Borel probability measure $\mu,$  the
integral (5.1) defines always a state-measure. A state-measure is
extremal if and only if it is a state-measure-morphism if and only
if $\mu = \delta_x$ for some point $x \in \Omega$, where
$\delta_x(M)=1$ iff $x \in M$ otherwise $\delta_x(M)=0,$ then
$m(f)=f(x).$

\end{ex}

We say that a net of  state-measures $\{m_\alpha\}$ {\it converges
weakly} to a  state-measure $m$ if $m(a) = \lim_\alpha m_\alpha(a)$
for every $a \in A.$

\begin{prop}\label{pr:6.4} The state spaces $\mathcal {SM}(A,u)$ and  $\mathcal {SMM}(A,u)$ are  compact
Hausdorff topological spaces.
\end{prop}

\begin{proof} If $\mathcal {SM}(A,u)$ is void, the statement is evident.
Thus suppose that $(A,u)$ admits at least one state-measure. For any
state-measure $m$ and any $x \in A$ we have by (2.1): $m(u\to x) =
m(u\vee_1 x \to x)=m(x)-m(u\vee_1 x).$  But $u\le u\vee_1 x,$ hence
$m(u\vee_1 x) \le m(u)=1$ so that $m(u\to x) \ge m(x)-1$ and
$$m(x) \le m(u\to x)+1.
$$
Therefore,
$$m(x) \le m(u\to x)+1 \le m(u\to^2 x)+2\le\cdots \le
m(u\to^{n-1}x)+n-1.
$$
Since $u$ is strong,  given $x \in A,$ let $n_x$ denote an integer
$n_x\ge 1$ such that $u\to^{n_x}x =1.$ Then $u\le (u\to^{n_x-1} x)$
and $m(u\to^{n_x-1} x)\le m(u)=1$. Consequently, $m(x)\le m(u\to
^{n_x-1} x) +n_x -1 \le n_x.$ Hence, $\mathcal {SM}(A,u) \subseteq
\prod_{x \in A}[0,n_x].$ By Tychonoff's Theorem, the product of
closed intervals is compact. The set of state-measures $\mathcal
{SMM}(A,u)$ can be expressed as an intersection of closed subsets of
$[0,\infty)^A$, namely of the following sets (for $x,y \in A$)

\begin{eqnarray*}
&M_{x,y}=\{m \in [0,\infty)^A: \ m(x\to y)=m(x\rightsquigarrow
y)=m(y)-m(x)\},\ x\le y,\\
& M_x=\{m\in [0,\infty)^A: m(x)\ge 0\},\quad \{m\in [0,\infty)^A:
m(u)=1\}.
\end{eqnarray*}
Therefore, $\mathcal {SM}(A,u)$ is a closed subset of the given
product of intervals, and hence, it is compact.

Similarly,  the set of state-measure-morphisms $\mathcal {SMM}(A,u)$
is a subset of $ \prod_{x \in A}[0,n_x]$ and it can be expressed as
an intersection of closed subsets of $[0,\infty)^A$, namely of the
following sets (for $x,y \in A$)

\begin{eqnarray*}
&M_{x,y}=\{m \in [0,\infty)^A: \ m(x\to y)=m(x\rightsquigarrow
y)=\max\{0,m(y)-m(x)\}\},\\
& M_x=\{m\in [0,\infty)^A: m(x)\ge 0\},\quad \{m\in [0,\infty)^A:
m(u)=1\}.
\end{eqnarray*}
Therefore, $\mathcal {SMM}(A,u)$ is a closed subset of the given
product of intervals, and hence, it is compact.
\end{proof}

\begin{prop}\label{pr:6.5}  Let $u$ be a strong unit of a pseudo-BCK
algebra $A$ and $m$ let be a measure on $A$.  Then $m$ vanishes on
$A$ if and only if $m(u)=0.$
\end{prop}

\begin{proof}  Assume $m(u) =0.$ Then $m(u\to x) =
m(u\vee_1 x \to x)=m(x)-m(u\vee_1 x).$  But $u\le u\vee_1 x,$ hence
$0\le m(u\vee_1 x) \le m(u)=0$ so that $m(x)= m(u\to x)$ and
$$m(x) = m(u\to x)=m(u\to^2 x)=\cdots=(m\to^n x)=m(1)=0$$
when $u\to^n x=1$ for some integer $n\ge 1.$

If now $m(u)>0,$ then $m$ does not vanish trivially on $A.$

\end{proof}

\begin{lemma}\label{le:6.7}  Let $m_1,m_2$ be  state-measure-morphisms on a
unital pseudo-BCK algebra $(A,u).$ If $\Ker_0(m_1)=\Ker_0(m_2),$
then $m_1=m_2.$

In addition, any state-measure-morphism cannot be expressed as a
convex combination of other state-measure-morphisms.
\end{lemma}

\begin{proof}  The set $m_1(A)=
\{m_1(a)\mid a\in A\}$ and $m_2(A)= \{m_2(a)\mid a\in A\}$ of real
numbers can be endowed with a total operation $*_\mathbb R$ such
$(m_1(A),*_\mathbb R,0)$ and $(m_2(A),*_\mathbb R,0)$ is a
subalgebra of the BCK algebra $([0,\infty),*_\mathbb R,0)$ in the
sense of \cite[Chap 5]{DvPu}, where $s*_\mathbb R t =
\max\{0,s-t\},$ $s,t \in [0,\infty).$ And the number $1$ is a strong
unit in all such algebras (for definition see \cite{DvPu}).

If we set $\hat m_1$ and $\hat m_2$ the state-measure-morphisms on
the quotient pseudo-BCK algebras $A/\Ker_0(m_1)$ and $A/\Ker_0(m_2)$
defined by $\hat m_i(a/\Ker_0(m_i))=m_i(a),$ we have again $\hat
m_i(A/\Ker_0(m_i))=m_i(A)$ for $i=1,2.$

Define a mapping $\phi:\ m_1(A) \to m_2(A)$ by $\phi(m_1(a)) =
m_2(a)$ ($a\in A$). It is possible to show that this is a BCK
algebra injective homomorphism. Due to \cite[Lem 6.1.22]{DvPu}, this
means that $m_1(A) = m_2(A)$ and $m_1(a)=m_2(a)$ for all $a\in A.$

Suppose now that $m= \lambda m_1 +(1-\lambda)m_2,$ where $m, m_1,
m_2$ are state-measure-morphisms and $\lambda \in (0,1).$ Then
$\Ker_0(m) \subseteq \Ker_0(m_1)\cap \Ker_0(m_2).$ Due to
Proposition \ref{pr:6.6}, all kernels
$\Ker_0(m),\Ker_0(m_1),\Ker_0(m_2)$ are maximal filters so that
$\Ker_0(m)=\Ker_0(m_1)=\Ker_0(m_2)$ and by the first part of the
present proof, $m=m_1 =m_2.$
\end{proof}

\begin{prop}\label{AD9}  Let $u$ be a strong unit of a pseudo-BCK
algebra $A$ and let  $J$ be a filter of $A,$ $J_0:=J\cap[u,1].$ Then
$J_0$ is a filter of the pseudo BCK-algebra
$([u,1],\le,\to,\squig,u,1).$ If $ F(J_0)$ is the filter of $A$
generated by $J_0,$  then
$$F(J_0) = F. \eqno(5.2) $$

Moreover, $J_0$ is maximal in $[u,1]$ if and only if so is $J$ in
$A.$
\end{prop}

\begin{proof}  Suppose that $J$ is a filter of $A.$ Then
$J_0 := J\cap [u,1]$ is evidently a filter of $[u,1].$

It is clear that  $F(J_0) \subseteq F.$

On the other hand, take $x \in J.$ Since $u$ is a strong unit, by
(2.3), there is an integer $n\ge 1$ such that $u\to^nx = 1 = u\to
(\cdots \to (u\to x)\cdots ).$  Set $x_n = u\vee_1 x$ and $x_{n-i}=
u\vee_1 (u\to^ix)$ for $i=1,\ldots,n-1.$ An easy calculus shows that
$x_i \in J_0$ for any $i=1,\ldots,n.$ Moreover, $u\to (u \to (\cdots
\to (u\to x)\cdots )) = x_1 \to (x_2\to( \cdots \to(x_n\to x)\cdots
))=1$ which by (2.2) proves $x \in F(J_0).$

Let now $J$ be a maximal filter of $A.$ Assume  that $F$ is a filter
of $[u,1]$  containing $J_0$ with $F \ne [u,1],$ and let $\hat F(F)$
be the filter of $A$ generated by $F.$ Then $F \subseteq \hat
F(F)\cap [u,1].$  If now $x \in \hat F(F)\cap [u,1],$ there are
$f_1,\ldots, f_n \in F$ such that $f_1\to(\cdots\to(f_n\to
x)\cdots)=1$ giving $x \in F.$ Hence, $F= \hat F(F)\cap [u,1].$

We assert that  $\hat F(F)$ is a  filter of $A$ containing $J,$ and
$\hat F(F)\ne A.$ If not, then $u \in \hat F(F)$ and therefore by
(2.2), there are $x_1,\ldots,x_n \in F$ such that $x_1\to(\cdots \to
(x_n\to u)\cdots)=1.$ If we set $z_n = x_n\vee_1 u$ and $z_{n-i}=
x_{n-i} \vee_1(x_i \to(\cdots \to(x_n\to u)\cdots )),$ for
$i=1,\ldots, n-1,$ then each $z_i$ belongs to $F$ and $z_1\to
(\cdots \to(z_n\to u)\cdots) =1$ which implies $u\in F$ that is a
contradiction.

The maximality of $J$ entails $J = \hat F(F).$  Since $J_0 \subseteq
F = \hat F(F) \cap [u,1] = J\cap [0,1]=J_0.$ That is, $J_0$ a
maximal filter of $[u,1]$ as it was claimed.

Assume now that $J_0$ is a maximal filter of $[u,1]$ and let $G\ne
A$ be a filter of $A$ containing $J.$ Then $G_0:=G\cap [u,1]$ is a
filter of $[u,1]$ containing $J_0,$ and by (5.2), $G=F(G_0).$ We
assert $u \notin G_0.$ Suppose the converse. Then $x \in G$ and for
any $x \in A$, there is an integer $n\ge 1$ such that $u\to^n x =
u\to (\cdots(u\to x)\cdots )=1$ proving $x \in G,$ so that
$A\subseteq G$ that is absurd.

The maximality of $J_0$ entails $J_0 = G_0$ and in view of (5.2), we
have $J=F(J_0)=F(G_0)=G,$ thus $J$ is a maximal filter of $A.$
\end{proof}

\begin{prop}\label{AD10}  Let $m$ be a state-measure  on a unital pseudo-BCK
algebra $(A,u), $ and let  $m_u$ be the restriction  of $m$ onto the
interval $[u,1].$ Then $m_u$ is a state-measure-morphism on
$([u,1],\le,\to,\squig).$

Let us have the following conditions:\\
$(a)$ $m$ is a state-measure-morphism on $(A,u).$\\
$(b)$ $m_u$  is a state-morphism on $[u,1].$\\
$(c)$ $\Ker_0(m)$ is a maximal filter of $(A,u).$\\
$(d)$ $\Ker_0(m_u)$ is a maximal filter of $[u,1].$

Then $(b), (c),(d)$ are mutually equivalent and $(a)$ implies each
of the conditions $(b),(c),$ and $(d).$

\end{prop}

\begin{proof}
Let $m_u$ be the restriction of $m$ onto $[u,1].$ Then $m_u$ is a
state-measure on $[u,1]$ and we define $\Ker_0(m_u)$ and
$\Ker_0(m)$.  Then  $\Ker_0(m_u)=\Ker_0(m)\cap [u,1]$ and due to
(5.2) we have

$$
F(\Ker_0(m_u))=\Ker_0(m).\eqno(5.3)
$$

$(a) \Rightarrow (b).$  It is evident.

$(b) \Leftrightarrow (d):$ It follows from  Theorem
\ref{psBCK-2.410}(d)-(e).

$(b)\Leftrightarrow (c):$  We have $\Ker_0(m_u)=\Ker_0(m)\cap
[u,1].$ Then by Proposition \ref{AD9}, we have the equivalence in
question.
\end{proof}

\begin{theo}\label{ADThm3}  Let $(A,u)$ be an unbounded pseudo-BCK
algebra that is downwards-directed and let $m$ be a state-measure on
$(A,u).$ Then there is a unique (up to isomorphism) abelian and
archimedean $\ell$-group $G$ with strong unit $u_G>0$ such that the
unbounded unital pseudo-BCK algebra $(A/\Ker_0(m),u/\Ker_0(m))$ is
isomorphic with the unbounded unital pseudo-BCK algebra
$(G^-,-u_G).$
\end{theo}

\begin{proof}  We define $m_u,$ the restriction of $m$ onto the
interval $[u,1].$  Due to Theorem \ref{AD6} and Theorem
\ref{psBCK-2.310}, the quotient $[u,1]/\Ker_0(m_u)$ can be converted
into an MV-algebra, and in view of $\Ker_0(m_u)=\Ker_0(m)\cap [u,1]$
we have that $[u,1]/\Ker_0(m_u)$ is isomorphic with $[u/
\Ker_0(m),1/\Ker_0(m)]=[u,1]/\Ker_0(m),$ so that both can be viewed
as isomorphic MV-algebras.  Let $G$ be an $\ell$-group guaranteed by
Theorem \ref{ADThm1} that is generated by $A/\Ker_0(m).$ Therefore,
$u_G:=-(u/\Ker_0(m))$ is a strong unit for $G$ and
$0_G:=1/\Ker_0(m)$ is the neutral element of $G.$ By the Mundici
famous theorem \cite{CDM}, the unital $\ell$-group $(G,u_G)$ is the
same for $[u,1]/\Ker_0(m_u)$ and $[u,1]/\Ker_0(m).$  If now $g \in
G^-,$ then $g = g_1+\cdots+g_n,$ where $g_1,\ldots, g_n\in
[u,1]/\Ker_0(m).$ The set of elements $g \in G^-$ such that $g \in
A/\Ker_0(m)$ is a pseudo-BCK algebra containing $A/\Ker_0(m).$ And
because $A/\Ker_0(m)$ generates $G,$ this implies that the
pseudo-BCK algebra $(G^-,-u_G)$ is isomorphic with the unital
pseudo-BCK algebra $(A/\Ker_0(m),u/\Ker_0(m)).$
\end{proof}

\begin{theo}\label{th:6.8} Let $m$ be a state-measure on a unital
pseudo-BCK algebra $(A,u)$ that is downwards-directed and let $m_u$
be the restriction of $m$ onto the pseudo-BCK algebra $[u,1].$ The
following
statements are equivalent:\\
$(a)$ $m$ is a state-measure-morphism on $(A,u).$\\
$(b)$ $m_u$  is a state-morphism on $[u,1].$\\
$(c)$ $\Ker_0(m)$ is a maximal filter of $(A,u).$\\
$(d)$ $\Ker_0(m_u)$ is a maximal filter of $[u,1].$\\
$(e)$ $m$ is an extremal state-measure on $(A,u).$\\
$(f)$ $m_u$ is an extremal state-measure on $[u,1].$
\end{theo}

\begin{proof} By Theorem \ref{AD10},  $(b), (c),(d)$ are mutually equivalent and
$(a)$ implies each of the conditions $(b),(c), (d).$  Theorem
\ref{ADThm2} entails that $(c)$ implies $(a).$  From Theorem
\ref{psBCK-2.410} we see that $(b)$ and $(f)$ are equivalent.
Proposition \ref{ADProp1} gives $(a)$ implies $(e).$

$(e) \Rightarrow (a).$  Let $m$ an extremal state-measure on
$(A,u).$ Define $\Ker_0(m),$ $A/\Ker_0(m),$ and $\hat
m(a/\Ker_0(m)):=m(a)$ $(a \in A).$  We assert that $\hat m$ is
extremal on the unital pseudo-BCK algebra $(A/\Ker_0(m),
u/\Ker_0(m)).$ Indeed, if $\hat m=\lambda \mu_1 +(1-\lambda)\mu_2,$
$0<\lambda <1,$ where $\mu_1$ and $\mu_2$ are two state-measures on
$(A/\Ker_0(m),u/\Ker_0(m)),$ then there are two state-measures $m_1,
m_2$ on $(A,u)$ such that $\hat m_1 =\mu_1$ and $\hat m_2 = \mu_2.$
Hence, $m=\lambda m_1 +(1-\lambda)m_2$ yielding $m_1 = m_2$ and
$\mu_1=\mu_2.$

Due to Theorem \ref{ADThm3}, $A/\Ker_0(m)$ is isomorphic with the
pseudo-BCK algebra $(G^-,u_G)$, where $G^-$ is the negative cone of
an abelian and archimedean $\ell$-group $G$ that is generated by
$A/\Ker_0(m)$ and the element $u_G:=-(u/\Ker_0(m))$ is a strong unit
for $G.$

Similarly as in the proof of Theorem \ref{ADThm2}, $\hat m$ can be
extended to a state-measure on $(G^-,u/\Ker_0(m))$ so that $s$ can
be extended to an additive function $s$ on the whole unital
$\ell$-group $(G,-(u/\Ker_0(m)))$ that is positive on $G^-$ and
$s(u/\Ker_0(m))=1.$  Moreover, $s$ is extremal on
$(G,-(u/\Ker_0(m)))$ which by \cite[Thm 12.18]{Goo} is possible if
and only if $\Ker_0(\hat m)=\{1/\Ker_0(m)\}$ is a maximal filter of
the unital pseudo-BCK algebra $(A/\Ker_0(m),u/\Ker_0(m)).$ Since the
mapping $a \mapsto a/\Ker_0(m)$ is surjective, we have that this
implies that $\Ker_0(m)$ is a maximal filter of $(A,u).$ By the
equivalence of $(c)$ and $(a)$ we have that $m$ is a
measure-morphism.
\end{proof}

As a direct consequence of Theorem \ref{th:6.8} and the
Krein--Mil'man theorem we have:

\begin{cor}\label{co:AD}  Let $(A,u)$ be a unital pseudo-BCK algebra
that is downwards-directed. Then
$$ \partial_e\mathcal{SM}(A,u) = \mathcal{SMM}(A,u) \eqno(5.4)$$
and every  state-measure on $(A,u)$ is a weak limit of  a net of
convex combinations of states-measure-morphisms.

\end{cor}

\section{Coherence, de Finetti Maps and Borel States}

In this section, we will generalize to pseudo-BCK algebras the
identity between de Finetti maps and Bosbach states, following the
results proved by K\"uhr and Mundici in \cite{Kumu} who showed that
de Finetti's coherence principle that has an origin in Dutch book
making, has a strong relationship with MV-states on MV-algebras.
Then we generalize this also for state-measures on unital pseudo-BCK
algebras that are downwards-directed.

Finally we present some open questions.

 We recall the following definition and notations used in
\cite{Kumu}. Let $A$ be a nonempty set, $[0,1]^A$ the set of all
functions $V:A \rightarrow [0,1]$ endowed with the product topology.
If $\mathcal {X}\subseteq [0,1]^A$, by $\rm conv \mathcal{X}$, $\rm
cl \mathcal{X}$ we denote the convex hull and respectively the
closure of $\mathcal{X}$. Also, if $\mathcal{X}$ is convex,
$\partial_e \mathcal{X}$ will denote the set of all extremal points
of $\mathcal{X}$.  We note that the weak topology of Bosbach states
is in fact the relativized product topology on $[0,1]^A.$

\begin{Def} {\rm (\cite{Kumu})}
Let $A'=\{a_1,a_2,\ldots,a_n\}$ be a finite subset of $A$. Then a
map $\beta:A'\longrightarrow [0,1]$ is said to be \emph {coherent}
over $A'$ if,
$$\mbox{for all}\ \sigma_1, \sigma_2,\ldots,\sigma_n\in
\mathbb R,\ \mbox{there is} \ V\in \mathcal{W}\ \mbox{s.t.}\
\sum_{i=1}^n \sigma_i(\beta(a_i)-V(a_i))\geq 0.\eqno(6.1)
$$
By a \emph {de
Finetti map} on $A$ we mean a function $\beta:A\longrightarrow
[0,1]$ which is coherent over every finite subset of $A$. We denote
by $\mathcal{F}_{\mathcal{W}}$ the set of all de Finetti maps on
$A$.
\end{Def}

An interpretation of (6.1) is as follows, \cite{Kumu}: Two players,
the bookmaker and the bettor, wager money on the possible occurrence
of elementary events $a_1,\ldots,a_n \in M.$ The bookmaker sets a
betting odd $\beta(a_i) \in [0,1],$ and the bettor chooses stakes
$\sigma_i \in \mathbb R.$ The bettor pays the bookmaker
$\sigma_i\beta(a_i),$ and will receive $\sigma_iV(a_i)$ from the
bookmaker's possible world $V.$ As scholars,  we can assume that
$\sigma_i$ could be also positive as well as negative. If the
orientation of money transfer  is given via bettor-to-bookmaker,
then (6.1) means that bookmaker's book should be {\it coherent} in
the sense that the bettor cannot choose stakes
$\sigma_1,\ldots,\sigma_n$ ensuring him to win money  in every $V
\in {\mathcal W}.$

Now let $A$  be a bounded pseudo-BCK algebra and denote by
$\mathcal{BS}(A)$ the set of Bosbach states on $A$ and by
$\mathcal{W}$ the set of state-morphisms on $A$. Note that,
according to Theorem \ref{psBCK-2.410}, $\mathcal{W}$ coincides with
the set of extremal Bosbach states, and due to the
Krein-Milman'theorem,

$$\mathcal{BS}(A) = \mbox{cl conv} \partial_e \mathcal{BS}(A)=
\mbox{cl conv} \mathcal{SM}(A).\eqno(6.2)
$$

\begin{theo}\label{Finet}  Let  $A$ be a bounded pseudo-BCK algebra and let
$\mathcal{W} = \mathcal{SM}(A) \ne \emptyset.$ Then
$$\mathcal{F}_{\mathcal{W}}=\mathcal{BS}(A).$$
\end{theo}

\begin{proof}
According to Theorem \ref{psBCK-2.410} and (3.1), $\mathcal {W}$ is
closed. Now we can apply  \cite[Prop 3.1]{Kumu} since we have
$\mathcal{W}\subseteq \mathcal{S}, \partial_e
{\mathcal{S}}=\mathcal{W}(\subseteq \mathcal{W})$, $\mathcal {W}$
closed. So we get $\mathcal{S}=\mathcal{F}_{\mathcal{W}}.$
\end{proof}

Theorem \ref{Finet} has an important consequence, namely that every
Bosbach state (if it exists) on a bounded pseudo-BCK algebra is a de
Finetti map coming from the set of $[0,1]$-valued functions on $A$
generated by the set of state-morphisms, and applying (6.2) we have
that this de Finetti maps is exactly the weak limit of a net of
convex combinations of state-morphisms.

There is also another relationship concerning the representability
of Bosbach states via integrals. We introduce the following notions,
see e.g. \cite[Sec 5]{Goo}.  Let  $\Omega$ be a nonempty compact
Hausdorff topological space. Let ${\mathcal B}(\Omega)$ be the Borel
$\sigma$-algebra of $\Omega$ generated by all open subsets of
$\Omega$ and any elements of ${\mathcal B}(\Omega)$ is said to be a
{\it Borel set}, and any $\sigma$-additive (signed) measure is said
to be a {\it Borel measure}.

Let  ${\mathcal P}(\Omega)$ denote all probability measures, that
is, all positive regular Borel measures $\mu \in {\mathcal
M}(\Omega)$ such that $\mu(\Omega) = 1.$ We recall that a Borel
measure $\mu$ is called regular if

$$\inf\{\mu(O):\ Y \subseteq O,\ O\ \mbox{open}\}=\mu(Y)
=\sup\{\mu(C):\ C \subseteq Y,\ C\ \mbox{closed}\}
$$
for any $Y \in {\mathcal B}(\Omega).$

Let now $A$ be a bounded pseudo-BCK algebra and let $\mathcal
W=\mathcal {SM}(A).$ Every element $a\in A$ determines a
(continuous) function $f_a: \ {\mathcal W} \to [0,1]$ via
$$ f_a(V) = V(a), \quad V \in {\mathcal W}.$$

We say that a mapping $s:\ A \to [0,1]$ is a {\it Borel state} (of
${\mathcal W}$) if there is a ($\sigma$-additive) probability
measure $\mu$ defined on the Borel $\sigma$-algebra of  the
topological space ${\mathcal W}$ generated by all open subsets of
${\mathcal W}$  such that
$$
s(a) =\int_{\mathcal W} f_a(V)\dx\mu(V).
$$
Let ${\mathcal B}_{\mathcal W}$ be the set of all Borel states of
${\mathcal W}.$

\begin{theo}\label{Finet1}  Let $A$ be
a bounded pseudo-BCK algebra. For any Bosbach state $s$ on $A$ there
is a Borel probability measure $\mu$ on $\mathcal B(\mathcal W)$
such that

$$
s(a) =\int_{\mathcal W} f_a(V)\dx\mu(V).\eqno(6.3)
$$

\end{theo}

\begin{proof}
Since ${\mathcal W}=\mathcal {SM}(A)$ is closed, see (3.1), by
\cite[Thm 4.2]{Kumu} we have ${\mathcal W} \subseteq {\mathcal
B}_{\mathcal W},$ $\partial_e {\mathcal B}_{\mathcal W} \subseteq
{\mathcal W}$ and ${\mathcal F}_{\mathcal W} = {\mathcal
B}_{\mathcal W}.$ Therefore, by Theorem \ref{Finet}, $\mathcal
{BS}(A) = {\mathcal B}_{\mathcal W},$ i.e. every Bosbach state is a
Borel state on $A.$
\end{proof}

We note that if we set $\Omega=\mathcal {BS}(A)$, then for any $a\in
A$, the function $\tilde a:\ \mathcal {BS}(A) \to [0,1]$ defined by
$\tilde  a(s)=s(a),$ $s\in \mathcal {SB}(A),$ is continuous.
Therefore, we can strength Theorem \ref{Finet1} as follows.

\begin{theo}\label{Finet2}  Let $A$ be
a bounded pseudo-BCK algebra. For any Bosbach state $s$ on $A$ there
is a unique Borel probability measure $\mu$ on $\mathcal B(\mathcal
{BS}(A))$ such that

$$
s(a) =\int_{\mathcal {SM}(A)} \tilde  a(x) \dx\mu(x).\eqno(6.4)
$$

\end{theo}

\begin{proof}  Suppose that the set of all Bosbach states on $A$ is
non empty.  Due to the Krein-Mil'man theorem, (6.2), the set of
extremal  Bosbach state is also nonempty and it coincides with the
set of state-morphisms. Denote by $F_0 := \bigcap\{\Ker(s):\ s \in
\mathcal {SM}(A)\}.$  In view of Propositions
\ref{psBCK-2.310}--\ref{psBCK-2.340}, $F_0$ is a normal ideal, and
similarly as in Theorem \ref{psBCK-2.310}, we can show that $A/F_0$
is an archimedean MV-algebra, and for any Bosbach state $s$ on $A$,
the mapping $\hat s(a/F_0)=s(a)$ $(a\in A)$ is an MV-state (=
Bosbach state) on  $A/F_0;$ we set $\bar a:= a/F_0$ ($a\in A$).
Moreover, the state spaces $\mathcal {BS}(A)$ and
$\mathcal{BS}(A/F_0)$ are affinely homeomorphic compact nonempty
Hausdorff topological spaces under the mapping $s\in \mathcal
{BS}(A)\mapsto \hat s\in \mathcal{SB}(A/F_0)$ (i.e. they are
homeomorphic in the weak topologies of states preserving convex
combinations of states). In addition, the compact subsets of
extremal Bosbach space are also homeomorphic under this mapping. Due
to \cite{Kro}, on the Borel $\sigma$-algebra $\mathcal B({\mathcal
{BS}(A)}),$ there is a unique Borel probability measure $\mu$ such
that
$$s(a)= {\hat s}(a/F_0) = \int_{\mathcal {SM}(A/F_0)} \tilde{\bar a} \dx\mu. \
$$
This integral can be rewritten identifying the compact spaces and
Borel $\sigma$-algebras into the form
$$ s(a) =\int_{\mathcal{SM}(A)} \tilde a(x)\dx\mu(x).$$
\end{proof}

It is interesting to note that de Finetti was a great propagator
only of probabilities as finitely additive measures. The result of
\cite{Kro} and formula (6.4) say that whenever $s$ is a Bosbach
state, it generates a $\sigma$-additive probability such that $s$ is
in fact an integral over this Borel probability measure. Thus
formula (6.4) joins de Finetti's  ``finitely additive probabilities"
with $\sigma$-additive measures on an appropriate Borel
$\sigma$-algebra.

We now generalize Theorem \ref{Finet} and Theorem \ref{Finet1} also
for unbounded pseudo-BCK algebras that are downwards-directed.

\begin{theo}\label{Finet6.1}  Let  $(A,u)$ be a  pseudo-BCK algebra that is
downwards-directed and let $\mathcal{W} = \mathcal{SMM}(A,u) \ne
\emptyset.$ Then
$$\mathcal{F}_{\mathcal{W}}=\mathcal{SM}(A,u).$$
\end{theo}

\begin{proof}
It follows from Theorem \ref{th:6.8} and using the same steps as
those in Theorem \ref{Finet}. \end{proof}

\begin{theo}\label{Finet6.2}  Let $(A,u)$ be
a  pseudo-BCK algebra that is downwards-directed. For any
state-measure $m$ on $(A,u),$ where $\mathcal{W} =
\mathcal{SMM}(A,u) \ne \emptyset,$ there is a Borel probability
measure $\mu$ on $\mathcal B(\mathcal W)$ such that

$$
m(a) =\int_{\mathcal W} f_a(V)\dx\mu(V).
$$

\end{theo}

\begin{proof} It follows from Theorem \ref{th:6.8} and it follows
the  analogous steps as those in Theorem \ref{Finet1}.
\end{proof}

\vspace{3mm}

During our study we have found that the following questions remained
open:

(1) Do arrows $\to/\Ker_0(m)$ and $ \rightsquigarrow/\Ker_0(m)$ on
$A/\Ker_0(m)$ coincide always  in Theorem \ref{psBCK-3.50} ?

(2) Is Theorem \ref{ADThm1} true without the assumption on
downwards-directness ?

(3) If $(A,u)$ is a unital pseudo-BCK algebra is then $A$
downwards-directed ?

\setlength{\parindent}{0pt}

\vspace*{3mm}
\begin{flushright}
\begin{minipage}{148mm}\sc\footnotesize
Lavinia Corina Ciungu\\
Polytechnical University of Bucharest
Splaiul Independen\c tei 113, Bucharest, Romania \&\\
State University of New York - Buffalo, 244 Mathematics Building,
Buffalo NY, 14260-2900, USA\\
{\it E--mail address}: {\tt
lavinia\underline{ }ciungu@math.pub.ro, lcciungu@buffalo.edu}\vspace*{3mm}

Anatolij Dvure\v censkij\\
Mathematical Institute, Slovak Academy of Sciences,
\v Stef\'anikova 49, SK-814~73 Bratislava, Slovakia\\
{\it E--mail address}: {\tt
dvurecen@mat.savba.sk}\vspace*{3mm}

\end{minipage}
\end{flushright}

\end{document}